\catcode`\@=11

\hsize=150mm \vsize=240mm
\hoffset=4mm \voffset=3mm
\pretolerance=500    \tolerance=1000  \brokenpenalty=5000

\catcode`\:=\active
\def:{\relax\ifhmode\ifdim\lastskip>\z@\unskip\fi\penalty\@M\ \fi\string:}

\catcode`\!=\active
\def!{\relax\ifhmode\ifdim\lastskip>\z@\unskip\fi\kern.2em\fi\string!}

\catcode`\?=\active
\def?{\relax\ifhmode\ifdim\lastskip>\z@\unskip\fi\kern.2em\fi\string?}

\def\^#1{\if#1i{\accent"5E\i}\else{\accent"5E #1}\fi}
\def\"#1{\if#1i{\accent"7F\i}\else{\accent"7F #1}\fi}

\def\Medbreak{\vskip-\lastskip\medbreak}
  
\let\bigf@nt=\tenrm     \let\smallf@nt=\sevenrm

\def\pc#1{\bigf@nt#1\smallf@nt}         \def\pd#1 {{\pc#1} }

\long\def\th#1 #2\enonce#3\endth{%
   \Medbreak
   {\pc#1} {#2\unskip} {\sl #3}\medskip}

\long\def\tha#1 #2\enonce#3\endth{%
   \Medbreak
   {\pc#1} {#2\unskip}\par\nobreak{\sl #3}\medskip}

\frenchspacing
\catcode`\@=12


\input amssym.tex

\def\cat{cat\'egorie}

\def\dm{\mathop{{\bf DM}^{-}}\nolimits}
\def\DM{\mathop{{\bf DM}}\nolimits}
\def\dme{\mathop{{\bf DM}^{-,\rm eff}}\nolimits}
\def\dmet{\mathop{{\bf DM}^{-}_{\acute et}}\nolimits}
\def\dman{\mathop{{\bf DM}^{-}_{an}}\nolimits}
\def\DMan{\mathop{{\bf DM}_{an}}\nolimits}
\def\dge{\mathop{{\bf DM}_{\rm gm}^{\rm eff}}\nolimits}
\def\dg{\mathop{{\bf DM}_{\rm gm}}\nolimits}
\def\sh{\mathop{\rm Shv}\nolimits}
\def\sp{\mathop{\rm Spec}\nolimits}
\def\ho{\mathop{\rm Hom}\nolimits}

\def\im{\mathop{\rm Im}\nolimits}
\def\gr{\mathop{\rm Gr}\nolimits}
\def\smc{\mathop{\rm Smcor}\nolimits}
\def\sm{\mathop{\rm Sm}\nolimits}
\def\sch{\mathop{\rm Sch}\nolimits}
\def\ke{\mathop{\rm Ker}\nolimits}
\def\Gal{\mathop{\rm Gal}\nolimits}

\def\pr{\mathop{\rm pr}\nolimits}

\def\om{{\bf\Omega}}
\def\D{{\bf D}}

\def\M{{\bf M}}

\def\C{{\bf C}}
\def\Q{{\bf Q}}
\def\R{{\bf R}}
\def\r{{\rm R}}
\def\Z{{\bf Z}}
\def\Ztr{{\bf Z}_{tr}}

\font\symboles=msam10 
\def\cqfd{{\hbox{\symboles\char'3}}}

\def\limind{\mathop{\oalign{\rm lim\cr
\hidewidth$\longrightarrow$\hidewidth\cr}}}

\def\build#1_#2^#3{\mathrel{\mathop{\kern 0pt#1}\limits_{#2}^{#3}}}

\def\hfl#1#2{\smash{\mathop{\hbox to 12mm {\rightarrowfill}}
\limits^{\scriptstyle#1}_{\scriptstyle#2}}}

\def\vfl#1#2{\llap{$\scriptstyle #1$}\left\downarrow
\vbox to 6mm{}\right.\rlap{$\scriptstyle #2$}}

\def\dfleche#1#2{\smash{\mathop{\hbox to
9mm{\rightarrowfill}}\limits^{\scriptstyle#1}_{\scriptstyle#2}}}

 \def\gfleche#1#2{\smash{\mathop{\hbox to
9mm{\leftarrowfill}}\limits^{\scriptstyle#1}_{\scriptstyle#2}}}

\def\diagram#1{\def\normalbaselines{\baselineskip=0pt
\lineskip=10pt\lineskiplimit=1pt} \matrix{#1}}

\def\relrightarrow{\mathrel{\hbox to 9mm{\rightarrowfill}}}

\newtoks\ref
\newtoks\AUTHOR
\newtoks\TITLE
\newtoks\BOOKTITLE
\newtoks\NOTE
\newtoks\PUBLISHER
\newtoks\YEAR
\newtoks\JOURNAL
\newtoks\VOLUME
\newtoks\PAGES
\newtoks\NUMBER
\newtoks\SERIES
\newtoks\ADDRESS
\def\pointir{\unskip . --- \ignorespaces}

\def\book{\leavevmode
{[\the\ref]\enspace}%
\the\AUTHOR\pointir
\the\TITLE, 
{\sl\the\SERIES}
{\bf \the\VOLUME}
\the\PUBLISHER,
\the\ADDRESS,
({\the\YEAR}).
\smallskip
\filbreak}

\def\article{\leavevmode
{[\the\ref]\enspace}%
\the\AUTHOR\pointir
\the\TITLE,
{\sl\the\JOURNAL},
{\bf\the\VOLUME},
({\the\YEAR}),
{no \the\NUMBER},
\the\PAGES.
\smallskip
\filbreak}

\def\articles{\leavevmode
{[\the\ref]\enspace}%
\the\AUTHOR\pointir
\the\TITLE,
{\sl\the\JOURNAL},
{\bf\the\VOLUME},
({\the\YEAR}),
\the\PAGES.
\smallskip
\filbreak}

\def\prep{\leavevmode
{[\the\ref]\enspace}%
\the\AUTHOR\pointir
\the\TITLE,
{\sl en pr\'eparation}.
\smallskip
\filbreak}

\def\preprint{\leavevmode
{[\the\ref]\enspace}%
\the\AUTHOR\pointir
\the\TITLE,
{\sl preprint}.
\smallskip
\filbreak}

\def\incollection{\leavevmode
{[\the\ref]\enspace}%
\the\AUTHOR \pointir
\the\TITLE,
{\bf in }\the\BOOKTITLE, 
{\sl \the\SERIES },
{\bf\the\VOLUME} ,
\the\PUBLISHER ,
\the\ADDRESS, 
({\the\YEAR}), 
\the\PAGES.
\smallskip
\filbreak}

\def\biblio#1{\vglue 15mm\centerline{\bf#1}\vskip 10mm}
\let\+=\tabalign
\def\signature#1\endsignature{\vskip 15mm minus 5mm\rightline{\vtop{#1}}}

\overfullrule=0pt
\parindent=0cm
\footline={\hfil\tenrm\folio\hfil}

\centerline {\bf REALISATION DE HODGE DES MOTIFS DE VOEVODSKY}
\bigskip\bigskip

\centerline {Florence LECOMTE et Nathalie WACH}

\medskip\centerline {IRMA - Strasbourg}
\bigskip\bigskip
{\leftskip 10pt \rightskip 10 pt {\sl
RESUME 

Pour un sous-corps du corps des complexes, nous d\'efinissons un foncteur de la cat\'egorie des motifs g\'eom\'etriques de Voevodsky vers la cat\'egorie des \Z-complexes de Hodge mixtes de Deligne [D74]. Les  filtrations par le poids et de Hodge sont repr\'esent\'ees par des foncteurs de troncature d'un complexe des poids \`a la Bondarko [Bo10] pour la premi\`ere et du complexe de De Rham [LW09] pour l'autre.
\par}}

\beginsection{Introduction}

\medskip

Apr\`es l'invention des motifs par Grothendieck, Deligne a d\'efini les structures de Hodge mixtes comme de tr\`es bonnes approximations de structures motiviques: elles forment des cat\'egories ab\'eliennes et tensorielles et sont construites \`a partir de structures \'el\'ementaires dites pures. Il est donc essentiel de construire les r\'ealisations de Hodge des motifs de Voevodsky, ce qu'a fait rationnellement Annette Huber en [H00] et [H04]. Dans cet article, qui fait suite \`a nos travaux sur les r\'ealisations de De Rham [LW09] et Betti [L08] nous construisons une r\'ealisation de Hodge \`a partir de foncteurs repr\'esentables.
\medskip
Soit $k$ un sous-corps du corps des nombres complexes. En [L08], l'une des auteurs a construit pour tout plongement $\sigma:k\hookrightarrow \C$ un foncteur de r\'ealisation topologique $t_{\sigma}$  de la cat\'egorie des complexes motiviques de Voevodsky $\dm(k)$ dans la cat\'egorie d\'eriv\'ee $D(Ab)$ des groupes ab\'eliens et d\'efini le foncteur de r\'ealisation de Betti
$${\bf H}_{\sigma}(\M,q)  = \ho_{D(Ab)}(t_{\sigma}(\M), (2i\pi)^q\Z).$$
En [LW09] nous avons construit un ind-motif de De Rham pour repr\'esenter le foncteur de r\'ealisation de De Rham dans la cat\'egorie des complexes motiviques non born\'es de Cisinski-D\'eglise $\DM(k)$ (cf [CD09])
$${\bf H}_{DR}(\M) = \ho_{\DM(k)}(\M,\om).$$
Notons que ces r\'ealisations sont d\'efinies comme cohomologies de complexes, \`a savoir 
$${\bf H}_{\sigma}^p(\M,q)  = H^p(\r^\cdot\ho_{D(Ab)}(t_{\sigma}(\M), (2i\pi)^q\Z)) \ \ {\rm et } \ \ {\bf H}_{DR}^p(\M) =H^p(\r^\cdot\ho_{\DM(k)}(\M,\om)).$$
Lorsque l'on se restreint \`a la cat\'egorie des motifs g\'eom\'etriques $\dg(k)$, le th\'eor\`eme de De Rham se g\'en\'eralise en un th\'eor\`eme de comparaison 
$${\bf H}_{DR}(\M)\otimes \C \simeq {\bf H}_{\sigma}(\M,q)\otimes \C$$
dont les fl\`eches sont \'egalement d\'efinies au niveau des complexes. 
\smallskip
Bondarko [Bo10] a muni la cat\'egorie $\dg(k)$ d'une structure de poids qui permet de d\'efinir sur  les foncteurs de r\'ealisations des filtrations par le poids $W$ compatibles \`a l'isomorphisme de comparaison.

Pour tout motif g\'eom\'etrique $\M$, choisissons un complexe des poids $\M^{(\cdot)}$. Alors le complexe diff\'erentiel gradu\'e
$$\oplus_{i\in\Z}\r^\cdot\ho_{\DM(k)} (\M^{(i)}, \om),$$ qui s'envoie dans le complexe $\r^\cdot\ho_{\DM(k)}(\M, \om)$,
calcule la r\'ealisation de De Rham de $\M$ et la filtration par le poids est d\'efinie par le fonteur de troncature \`a droite sur la premi\`ere variable du complexe de motifs $\M^{(\cdot)}$
$$W_n \r^\cdot\ho_{\DM(k)}(\M,\om )= \im\left(\oplus_{i\leq n}\r^\cdot\ho_{\DM(k)} (\M^{(i)}, \om))\rightarrow \r^\cdot\ho_{\DM(k)}(\M, \om)\right).$$
On a de m\^eme, avec les notations \'evidentes 
$$W_n\r^\cdot\ho_{D(Ab)}(t_{\sigma}(\M), (2i\pi)^q\Z)= \im \left(\oplus_{i\leq n}\r^\cdot\ho_{D(Ab)} (t_{\sigma}(\M^{({i})}), (2i\pi)^q\ )\rightarrow \r^\cdot\ho_{D(Ab)}(t_{\sigma}(\M), (2i\pi)^q\Z)\right).$$

Sur $\r^\cdot\ho(\M^{(\cdot)}, \om)$ on tronque la deuxi\`eme variable en consid\'erant la filtration b\^ete de $\om$ pour d\'efinir une filtration $F$, dite filtration de Hodge. Celle-ci s'\'etend \`a
$\r^\cdot\ho(\M^{(\cdot)}, \om)\otimes\C$ o\`u elle induit une filtration oppos\'ee $\bar F$. Notre principal r\'esultat est 

\th THEOREME 0.1
\enonce Pour tout motif g\'eom\'etrique $\M$ d\'ecompos\'e en un complexe des poids $\M^{(\cdot)}$, les donn\'ees
$$\left(\r^\cdot\ho_{D(Ab)}(t_{\sigma}(\M), (2i\pi)^q\Z), \ \r^\cdot\ho(\M, \om),\ \r^\cdot\ho(\M, \om)\otimes\C\right)$$
munies des filtrations induites par les images des foncteurs de troncatures sur $\M^{(\cdot)}$ et $\om$ d\'efinissent un foncteur vers la cat\'egorie des $\Z$-complexes de Hodge mixtes au sens de Deligne [D74].
\endth
Le th\'eor\`eme signifie, entre autres, que les cohomologies des complexes $\r^\cdot\ho$ ne d\'ependent pas du choix du complexe des poids et que la filtration de Hodge est bien d\'efinie, bien que, a priori, en tronquant b\^etement $\om$ on perde l'invariance d'homotopie et on sorte de la cat\'egorie des complexes motiviques.
\smallskip
Comme chez Deligne [D74] il est n\'ecessaire de d\'efinir les r\'ealisations au niveau des complexes $\r^\cdot\ho$, afin d'obtenir une filtration de Hodge sur ${\bf H}_{DR}(\M)$ et non sur son gradu\'e par les poids. N\'eanmoins, les d\'emonstrations consistant \`a v\'erifier des propri\'et\'es cohomologiques, les r\'esultats seront \'enonc\'es en termes de foncteurs cohomologiques.

Bien qu'en [L08] et [LW09], les foncteurs de r\'ealisations aient \'et\'e d\'efinis sur la cat\'egorie $\dm(k)$ des complexes motiviques de Voevodsky, pour le th\'eor\`eme de comparaison ou la d\'eg\'en\'erescence des suites spectrales associ\'ees aux filtrations, nous avons besoin de r\'esultats de finitude et devons donc nous restreindre \`a la cat\'egorie des motifs g\'eom\'etriques $\dg(k)$ introduite par Voevodsky. 

Dans la cat\'egorie $\dm(k)$,
la cat\'egorie $\dg(k)$ est la sous-cat\'egorie 
pleine additive \'epaisse engendr\'ee par les motifs $\M(X)$ des sch\'emas 
projectifs lisses [MVW 14.1]. Cela signifie que la cat\'egorie $\dg(k)$ est 
construite \`a partir des sommes finies de motifs de sch\'emas projectifs 
lisses avec les deux propri\'et\'es suivantes:

    {\leftskip1cm
\itemitem{-}
    les facteurs directs des motifs g\'eom\'etriques sont g\'eom\'etriques;
  \itemitem{-}   pour tout triangle distingu\'e 
  $A\rightarrow B\rightarrow C\rightarrow A[1]$, si deux des trois motifs $A$, 
  $B$ et $C$ sont g\'eom\'etriques, le troisi\`eme l'est aussi.
 \par}

\tha PRINCIPE 0.2.
  \enonce Soit $H: \dme(k) \rightarrow {\cal A}$ un foncteur cohomologique
  vers une cat\'egorie ab\'elienne ${\cal A}$. Soit ${\cal B}$ une 
  sous-cat\'egorie pleine de $\cal A$, ab\'elienne, stable par facteur direct 
  et extension.
  Si $H(\M (X)[n])$ est un objet de ${\cal B}$ pour tout sch\'ema $X$ lisse et 
  projectif sur  $k$ et tout entier $n$, alors $H$ induit un foncteur 
  $$H : \dge(k) \rightarrow {\cal B}.$$
  Si de plus la cat\'egorie ${\cal B}$ est tensorielle, si $H$ est multiplicatif 
  et 
  $H({\bf Z}(1))$ est inversible dans ${\cal B}$, alors $H$ induit un foncteur
  $$H: {\bf DM}_{\rm gm}(k) \rightarrow {\cal B}.$$
  \endth

Dans cet \'enonc\'e, l'hypoth\`ese $\cal A$ ab\'elienne est n\'ecessaire pour consid\'erer un foncteur cohomologique.
\medskip

Nous commen\c cons par rappeler  les r\'esultats de Voevodsky et Bondarko. Puis nous revoyons bri\`evement la construction du motif de De Rham et d\'efinissons les filtrations par le poids et de Hodge. Ensuite nous construisons la r\'ealisation de Betti qui a \'et\'e r\'esum\'ee en [L08]. Elle est \'equivalente \`a celle construite par Ayoub en th\'eorie homotopique des sch\'emas [A10]. Finalement nous comparons ces deux r\'ealisations pour construire la r\'ealisation de Hodge, qui nous permet d'\'etendre aux motifs g\'eom\'etriques les cohomologies de Deligne-Beilinson. Nous terminons en montrant que rationnellement nos r\'ealisations co\"incident avec celles pr\'ec\'edemment construites par Huber [H00]. 

\smallskip
{\pc REMERCIEMENTS}: ce travail a b\'en\'efici\'e de discussions avec J. Wildeshaus sur les travaux de Bondarko et avec A. Huber sur la filtration de Hodge. Qu'ils en soient remerci\'es.

\th{CONVENTIONS}
\enonce Par la suite tous les corps sont suppos\'es de caract\'eristique $0$ et 
les sch\'emas sont s\'epar\'es de type fini sur un corps. On note $\sm(k)$ la \cat\ des sch\'emas lisses sur k, dont les morphismes sont les morphismes de sch\'emas. 
\endth

\beginsection{1. Motifs  de Voevodsky }

La cat\'egorie motivique dans laquelle nous travaillons est la cat\'egorie triangul\'ee
$\dm(k)$ de Voevodsky  ([V-TCM]) dont nous appelons les objets complexes motiviques. 
Cette cat\'egorie  est
obtenue par une s\'erie de  localisations \`a partir de la cat\'egorie des complexes de faisceaux sur le site $\smc(k)$ des correspondances finies et topologie de Nisnevich.

\beginsection{1.1. Rappels sur les sites et topos} 

Rappelons quelques notions utiles de [SGA4]. 
 Un site est une cat\'egorie 
munie d'une topologie de Grothendieck [SGA4 II 1.15]. Pour un site $\cal C$, on 
note $\widehat{\cal C} $ (resp. $\widetilde{\cal C}$) la cat\'egorie des 
pr\'efaisceaux (resp. faisceaux) d'ensembles du site $\cal C$. La cat\'egorie 
$\widetilde{\cal C}$ est appel\'ee topos associ\'e au site $\cal C$. Les 
cat\'egories $\widehat{\cal C} $ et $\widetilde{\cal C}$ sont munies de 
foncteurs canoniques $\cal C \rightarrow \widehat{\cal C}$ et
$\cal C \rightarrow \widetilde{\cal C}$ qui \`a un objet $C$ de  $\cal C$ 
associe respectivement le pr\'efaisceau et le faisceau repr\'esent\'es par 
l'objet $C$. Si $\cal C$ et $\cal C'$ sont deux sites, un foncteur 
$u: {\cal C}\rightarrow {\cal C'}$ entre les cat\'egories sous-jacentes induit 
par composition un foncteur 
$\hat u^* : \widehat {\cal C'}\rightarrow\widehat{\cal C} $. 
Ce foncteur $\hat u^*$ admet un adjoint \`a gauche 
$u_! :\widehat  {\cal C}\rightarrow\widehat{\cal C'}$ qui prolonge le foncteur 
d'origine $u$ 
dans le sens que le diagramme suivant commute
$$\diagram{
&{\cal C} & \hfl{u}{} & {\cal C}' \cr
&\vfl{}{} & & \vfl{}{}\cr
&\widehat{\cal C}& \hfl{u_!}{}& \widehat{\cal C}' \cr
}\leqno(1.1.1)$$
o\`u les fl\`eches verticales sont les foncteurs canoniques.

Le foncteur de prolongement $u_!$ est d\'efini de la fa\c con suivante: pour 
tout objet $C'$ 
de ${\cal C}'$ on note $I^{C'}_u$ la cat\'egorie des couples $(S,f)$ o\`u $S$ 
est un objet de $\cal C$ et $f$ un morphisme $f: C'\rightarrow u(S)$ dans
$\cal C'$ [loc. cit. I 5] et$(I^{C'}_u)^{\rm op}$ la cat\'egorie oppos\'ee.
On pose  
$$u_! F: C'\mapsto \limind_{(I^{C'}_u) ^{\rm op}}F\circ {\rm pr}_{C'} ( \ )$$ 
o\`u $\pr_{C'}$ est le foncteur de $I^{C'}_u$ vers $\cal C$ qui au couple 
$(S,f)$ 
associe l'objet $S$.\goodbreak

Si de plus, le foncteur $u$ est continu, c'est-\`a-dire que pour tout
 faisceau $G$ sur $\cal C'$ le pr\'efaisceau $C \mapsto G\circ u(C)$ est 
 un faisceau sur $\cal C$, alors le foncteur $u^*$ induit un foncteur 
 $u_s:\widetilde {\cal C'}\rightarrow\widetilde{\cal C} $ et ce foncteur $u_s$
 admet un adjoint \`a gauche $u^s$ qui prolonge $u$ 
 
$$\diagram{
&{\cal C} & \hfl{u}{} & {\cal C}' \cr
&\vfl{}{} & & \vfl{}{}\cr
&\widetilde{\cal C}& \hfl{u^s}{}& \widetilde{\cal C}' \cr
}\leqno (1.1.2)$$
Le foncteur $u^s$ est le compos\'e du foncteur de prolongement $u_!$ d\'efini 
plus haut avec le foncteur 
faisceau associ\'e. Par construction, il est exact \`a droite et commute aux
limites inductives. S'il est de plus exact \`a gauche, il est foncteur image 
inverse $u^s=\Phi^*$
d'un morphisme de topos $\Phi : \widetilde {\cal C'}\rightarrow\widetilde{\cal 
C} $
(loc. cit. IV.3.1.). Suivant toujours [SGA4], on note $\widetilde{\cal C}_{Ab}$ 
le topos ab\'elien associ\'e au site $\cal C$, c'est-\`a-dire la cat\'egorie 
des  faisceaux en groupes ab\'eliens sur $\cal C$.
\smallskip
Rappelons qu'un foncteur d'une cat\'egorie triangul\'ee $\cal T$ vers une cat\'egorie 
ab\'elienne $\cal A$ est dit cohomologique [Ver77] s'il transforme tout triangle 
distingu\'e en suite exacte. En particulier, pour un foncteur cohomologique $H$ et $p$ un entier, nous notons $H^p$ 
 le foncteur de $\cal T$ vers $\cal A$ 
d\'efini par $H^p(M) = H(M[-p])$;
 ceci  permet 
d'associer \`a
tout triangle distingu\'e une suite exacte longue.
Les foncteurs Hom sont des foncteurs cohomologiques. 

\beginsection{1.2. Les correspondances finies}

 Le groupe $Cor(X,Y)$ des  
correspondances finies entre deux sch\'emas lisses $X$ et $Y$ est le groupe ab\'elien libre engendr\'e par les sous-vari\'et\'es  
 ferm\'ees irr\'eductibles de $X \times_{\sp (k)} Y$ qui sont finies et surjectives 
 sur une composante irr\'eductible de $X$. Cette d\'efinition reste valable pour un sch\'ema $Y$ quelconque.

Les correspondances finies se comportent mieux que les cycles classiques: il existe des morphismes image inverse et image directe pour tous les morphismes entre sch\'emas lisses et elles se composent comme les correspondances de Grothendieck [Ma68]. Elles permettent de d\'efinir la \cat\ $\smc(k)$ des correspondances finies, dont les objets sont les sch\'emas lisses sur $k$ et les morphismes, les correspondances finies. La \cat\ $\smc(k)$ est additive, pour l'union disjointe, et tensorielle, pour le produit fibr\'e sur $\sp k$. Le foncteur canonique $\sm(k) \rightarrow \smc(k)$ qui envoie tout morphisme sur son graphe est compatible \` a ces structures. On appelle pr\'efaisceau avec transferts un 
foncteur contravariant de $\smc(k)$
 vers la cat\'egorie ${\it Ab}$ des groupes ab\'eliens; un faisceau de Nisnevich
 avec transferts est un pr\'efaisceau avec transferts qui est un faisceau pour 
 la topologie de Nisnevich, la topologie totalement d\'ecompos\'ee de [N89],
 interm\'ediaire entre la topologie de Zariski et la topologie \'etale.

La topologie de Nisnevich  munit les cat\'egories $\sm(k)$ et $\smc(k)$ ([BV08]4.3.) de topologie de Grothendieck. Le fait que $\smc(k)$ soit un site pour la topologie de Nisnevich provient de ce que l'image inverse d'un point (anneau hens\'elien) par une correspondance finie est un point. Suivant Voevodsky, on pr\'ef\'erera  noter $\sh_{Nis}(\smc(k))$ et $\sh_{\acute et}(\smc(k))$ les topos ab\'eliens des correspondances finies. Le foncteur canonique $ \sm(k) \rightarrow \smc (k)$ est continu (un pr\'efaisceau sur $\smc(k)$ est un faisceau si c'est un faisceau sur $\sm(k)$) et cocontinu car les cribles sur $\sm(k)$ et $\smc(k)$ sont les m\^emes.
  
\beginsection{1.3. Les cat\'egories motiviques} 

1.3.1. {\sl Faisceaux de Nisnevich avec transferts}

Les cat\'egories motiviques sont construites \`a partir de la \cat\ 
$\sh_{Nis}(\smc(k))$ des faisceaux de Nisnevich avec transferts.
Le faisceau  $\Ztr (X)$ (not\'e $L(X)$ dans 
 [V-TCM]) est le faisceau de Nisnevich repr\'esent\'e par le sch\'ema $X$ sur $\smc(k)$:
pour tout sch\'ema lisse $U$, on a $\Ztr(X)(U) = Cor(U,X)$. Notons que le faisceau 
$\Ztr (X)$ est d\'efini pour X quelconque.
Le produit des sch\'emas permet de d\'efinir le produit tensoriel des faisceaux avec transferts 
$$\Ztr(X)\otimes \Ztr(Y) = \Ztr(X \times_{\sp (k)} Y)$$
 pour toute paire de sch\'emas lisses $(X,Y)$.
Soulignons une propri\'et\'e   de la topologie de Nisnevich ([V-TCM] 
 Prop 3.1.3):

 \smallskip
\th PROPOSITION 1.3.1.1. {\rm (Voevodsky [V-TCM])}
\enonce Soit $X$ un sch\'ema lisse sur $k$ et 
${\cal U} = \lbrace U_i \rightarrow X \rbrace $ 
un recouvrement de Nisnevich de $X$. Notons $U$ l'union disjointe 
$U = \coprod U_i$ et $\check N ({\cal U}/X)$ le complexe de faisceaux 
$$\cdots \rightarrow \Ztr (U\times_X U) \rightarrow \Ztr (U) 
\rightarrow \Ztr (X) \rightarrow 0$$
avec les diff\'erentielles \'egales \`a la somme altern\'ee des morphismes 
induits par les projections.

Alors le complexe $\check N ({\cal U}/X)$ est acyclique pour la topologie 
de Nisnevich.
\endth

Comme tout sch\'ema lisse de type fini sur un corps 
 peut \^etre recouvert par une famille de sch\'emas lisses quasi-projectifs,
 cette proposition permet de r\'esoudre  les faisceaux $\Ztr (X)$, pour $X$ sch\'ema
 lisse de type fini par un complexe form\'e de sommes
  $\coprod_{\alpha} \Ztr (X_{\alpha})$ o\`u les sch\'emas $(X_{\alpha})$ sont 
  quasi-projectifs lisses. 
  C'est pourquoi dans nos constructions nous pourrons supposer que les sch\'emas
  sont quasi-projectifs.

La proposition (1.3.1.1.) reste valable en topologie \'etale mais pas en topologie de 
Zariski (loc. cit.). 

\smallskip\goodbreak
1.3.2. {\sl ${\bf A}^1$- localisation.}

La cat\'egorie $\dme(k)$ des complexes motiviques effectifs est la localisation de la cat\'egorie d\'eriv\'ee 
 $D^-(\sh_{Nis}(\smc(k))$ des complexes, born\'es sup\'erieurement, de 
 faisceaux de Nisnevich avec transferts par la sous-cat\'egorie
 \'epaisse engendr\'ee par les complexes du type
 \hbox{$\Ztr(X\times_{\sp (k)} {\bf A}^1_k) \rightarrow \Ztr(X)$.}
  On appelle {\sl $ {\bf A}^1_k$-\'equivalence} tout morphisme de 
$D^-(\sh_{Nis}(\smc(k)))$ qui induit un isomorphisme sur $\dm(k)$.
On note
  $\M(X)$  le complexe motivique  associ\'e au sch\'ema lisse $X$: il  est repr\'esent\'e dans $\dme(k)$ par le complexe 
singulier simplicial $ C_*(\Ztr(X))$ associ\'e \`a $X$, aussi appel\'e {\sl complexe 
de Suslin} du sch\'ema $X$. Pour un faisceau  $F$, le complexe $C_*(F)$ est le 
complexe de faisceaux d\'efini par 
 $$C_n(F)(X) = F(X\times \Delta^n)$$
  o\`u $\Delta^\cdot$ est le sch\'ema cosimplicial standard 
$\Delta^n = \sp k[z_0, \dots, z_n]/(\sum_{0\leq i \leq n} z_i - 1)$ et la 
diff\'erentielle est induite par la somme altern\'ee des morphismes de coface.

Le produit sur $\smc(k)$ se transporte sur $\dme(k)$ et on a pour toute paire $(X,Y)$ de sch\'emas lisses 
$$\M(X) \otimes \M(Y) = \M ( X \times_{\sp (k)} Y).$$
Le motif $\M({\bf P}^1)$ de la droite  projective  se scinde en 
$\M({\bf P}^1) = \Z \oplus \Z(1)[2]$ o\`u $\Z=\M(\sp k)$  est le motif du point et 
$\Z(1)$ est le motif de Tate, motif r\'eduit de ${\bf G}_m$. La \cat\ $\dm(k)$ des complexes motiviques est obtenue \` a partir de la \cat\ $\dme(k)$ en inversant le motif de Tate. Le th\'eor\`eme de simplification  ("cancellation theorem") de Voevodsky (cf [MVW] 16.25) permet 
d'identifier la cat\'egorie des complexes motiviques effectifs \` a une sous-cat\'egorie pleine de $\dm(k)$. Pour tout entier $q$ et tout complexe motivique $\M$ on note 
$\M(q)$ le produit $\M \otimes \Z (q)$.

\smallskip
{\pc REMARQUES}: hors le th\'eor\`eme de simplification, ces constructions sont formelles et restent valables en topologie \'etale sur $\smc(k)$ mais aussi sur $\smc(\C)$ pour la topologie analytique. Le th\'eor\`eme de simplification  est trivial en topologie analytique: le motif de Tate analytique est le complexe concentr\'e en degr\'e $0$ , $\Z(q)_{an}= (2i\pi)^q$ (cf ci-dessous, proposition 3.1.7.). 

Outre la cat\'egorie $\dmet(k)$ (loc. cit.) on peut construire la cat\'egorie triangul\'ee  $\dman$ \`a partir de la cat\'egorie $\smc(\C)$ et de la topologie analytique. Elle est \'equip\'ee d'un foncteur induit par le changement de topologie 
$$t_{an}: \dm(\C) \rightarrow \dman. \leqno (1.3.2.1)$$

Sp\'ecifique \`a la topologie de Nisnevich est le th\'eor\`eme d'invariance d'homotopie qui permet d'identifier les objets ${\bf A}^1$-locaux.
\smallskip
1.3.3. {\sl Complexes ${\bf A}^1$-locaux en topologie de Nisnevich}

On dit qu'un (pr\'e)-faisceau $F$ est {\sl invariant par homotopie} si pour tout sch\'ema 
lisse $X$, la projection $X \times_{\sp (k)}{\bf A}^1_k\rightarrow X$ induit un 
isomorphisme $F(X) \simeq F(X \times_{\sp (k)}{\bf A}^1_k)$.
Un r\'esulat fondamental de Voevodsky est le th\'eor\`eme d'invariance d'homotopie:

\th THEOREME 1.3.3.1. {\rm (Voevodsky [V-TCM] 3.1.12)}
\enonce Soit $F$ un faisceau de Nisnevich invariant par homotopie. Alors le faisceau de cohomologie associ\'e est \'egalement invariant par homotopie et on a pour tout sch\'ema 
$X$ lisse
et tout entier $i$ des isomorphismes 
$$\matrix{H^i_{Zar} (X,F) &\simeq &H^i_{Zar} (X \times_{\sp (k)}{\bf A}^1_k,F)\cr
\mid\!\wr&&\mid\!\wr\cr
 H^i_{Nis} (X,F) & \simeq &H^i_{Nis} (X \times_{\sp (k)}{\bf A}^1_k,F).}$$
\endth
\smallskip
Le th\'eor\`eme d'invariance d'homotopie
(1.3.3.1.) permet d'identifier $\dme(k)$ 
 \` a une sous-cat\'egorie pleine de la \cat\ d\'eriv\'ee $D^-(\sh_{Nis}(\smc(k)))$: en effet,
$\dme(k)$ en est la sous-cat\'egorie form\'ee des  complexes, born\'es sup\'erieurement, de faisceaux de Nisnevich avec transferts, qui sont \` a cohomologie invariante par homotopie.

Par la suite nous consid\'erons des complexes de faisceaux avec transferts 
 $L^\cdot$ qui sont ${\bf A}^1$-locaux, c'est-\`a-dire tels que pour tout 
 complexe motivique on a
 $$\ho_{D^-(\sh_{Nis} (\smc (k)))}(\M, L^\cdot)= \ho _{\dme(k)}(\M, L^\cdot).$$
En travaillant dans la cat\'egorie $D^-(\sh_{Nis}(\smc(k)))$ nous nous ramenons 
\`a d\'eriver des foncteurs de la cat\'egorie ab\'elienne des faisceaux de 
Nisnevich avec 
transferts. Cette cat\'egorie a assez d'injectifs [MVW 6.19] et les foncteurs 
$\rm Ext$ sont les d\'eriv\'es des foncteurs $\ho$. Plus pr\'ecis\'ement, nous 
notons
$\r^\cdot\ho$ le bifoncteur  d\'eriv\'e 
$$ D^-(\sh_{Nis} (\smc (k))) \times D^+(\sh_{Nis}(\smc(k))) \rightarrow D({\cal A}b)$$
d\'efini par 
$\r^\cdot\ho(M,N) = \ho^\cdot (M,I)$ 
o\`u  $M$ (resp. $N$) est un complexe de faisceaux de Nisnevich avec transferts 
born\'e sup\'erieurement (resp. inf\'erieurement), $I$ est une r\'esolution 
injective de $N$ et $\ho^\cdot(-,-)$ est le complexe 
$n \mapsto \ho(-,-[n])$. Le foncteur $\rm Ext$ est le 
foncteur cohomologique associ\'e et on a pour tout $M$ de $D^-(\sh_{Nis} \smc k)$ et
 tout complexe born\'e $N$ 
$${\rm Ext}^i(M,N) \simeq H^i(\r^\cdot\ho(M, N)) \simeq \ho_
{D^-(\sh_{Nis}(\smc k))}(M,N[i]).$$
Nous utilisons abondamment le r\'esultat suivant:

\smallskip
\th THEOREME 1.3.3.2. {\rm (Voevodsky [V-TCM])}
\enonce Si $L^\cdot$ est un complexe born\'e ${\bf A}^1 $-local de $\dme(k)$,
 alors 
pour tout motif $\M(X)$ d'un sch\'ema $X$ lisse sur $k$, on a des isomorphismes
$$\matrix{
\ho_{\dm(k)}(\M(X), L^\cdot[i]) &\simeq &\ho_{D^-(\sh_{Nis} (\smc (k)))}(\M(X), 
L^\cdot[i]) &\simeq&  {\Bbb H}^i_{Nis}(X,L^\cdot)\cr
&\simeq &\ho_{D^-(\sh_{Zar} (\smc (k)))}(\M(X), 
L^\cdot[i]) &\simeq& {\Bbb H}^i_{Zar}(X,L^\cdot)\cr
}$$
o\`u ${\Bbb H}^i_{Nis}$ (resp.${\Bbb H}^i_{Zar}$) d\'esigne l'hypercohomologie de 
Nisnevich (resp. Zariski) des complexes de faisceaux.
\endth

Le lemme ci-dessous permet de passer de la  cat\'egorie 
$D^-(\sh_{Nis} (\smc (k)))$ \`a $\dme(k)$.

\th LEMME d'homotopie 1.3.3.2.
\enonce Soit $F: D^-(\sh_{Nis} (\smc (k)))\rightarrow A$ ( resp. $ F: D^-(\sh_{Nis} (\smc (k)))\rightarrow T$) un foncteur cohomologique dans une cat\'egorie ab\'elienne $A$ (resp. foncteur exact dans une cat\'egorie triangul\'ee $T$) qui v\'erifie la propri\'et\'e suivante:

pour tout sch\'ema $X$ lisse sur $k$, la premi\`ere projection $\pi_X: X\times_k {\bf A}^1 
\rightarrow X$ induit un isomorphisme
$$F(\Z_{tr}(X)) \simeq F(\Z_{tr}(X\times_k {\bf A}^1)).$$
Alors $F$ se factorise en un foncteur cohomologique $F:\dme(k) \rightarrow A$ (resp. foncteur exact  $F:\dme(k) \rightarrow T$).
\endth
{\pc DEMONSTRATION}: Il suffit de montrer que si $f: K\rightarrow K'$ est une ${\bf A}^1$-\'equivalence de complexes de faisceaux avec transferts, alors l'image $F(Cf)$ du c\^one de $f$ est nulle. Par d\'efinition, le c\^one $Cf$ est dans la plus petite cat\'egorie \'epaisse contenant le c\^one $C\pi_X$ de $\pi_X$ et stable par somme directe.  Etant cohomologique (resp. exact), le foncteur $F$  commute aux sommes finies et facteurs directs et on se ram\`ene \`a montrer que l'image $F(C\pi_X)$ est nulle, ce qu'implique l'hypoth\`ese.
\hfill\cqfd

\smallskip
1.3.4. {\sl Motifs g\'eom\'etriques}

La cat\'egorie $\dge(k)$ des motifs g\'eom\'etriques effectifs est la sous-\cat\ \'epaisse de $\dme(k)$ engendr\'ee par les motifs $\M(X)$ des sch\'emas lisses. Comme le corps $k$ v\'erifie la r\'esolution des singularit\'es, la cat\'egorie $\dge(k)$ est engendr\'ee par les motifs des sch\'emas projectifs et lisses et elle contient les motifs de tous les sch\'emas sur $k$. Cette \cat\ est \'egalement  construite par double localisation (invariance d'homotopie et Mayer-Vietoris) de la cat\'egorie homotopique des complexes born\'es de $\smc(k)$. C'est cette deuxi\`eme construction qu'utilise Bondarko [Bo09] pour munir $\dge(k)$ d'une structure diff\'erentielle gradu\'ee. La cat\'egorie des motifs g\'eom\'etriques $\dg (k)$ est obtenue par inversion du motif de Tate. 

Dans la \cat\ des motifs g\'eom\'etriques $\dg(k)$, Voevodsky d\'efinit un Hom interne 
$\underline{\rm Hom}$ et une dualit\'e $\M^* = \underline{\rm Hom}(\M,\Z)$. Il associe \'egalement \` a tout sch\'ema $X$ un motif \`a support compact $\M^c (X)$, qui v\'erifie 
$\M^c (X)=\M(X)$, si $X$ est un sch\'ema projectif. On a, pour tout sch\'ema $X$ lisse de dimension $n$, la relation
([V-TCM] 4.3.2.)
$$ \M(X)^* = \M^c(X)(-n)[-2n].²$$

\smallskip
1.3.5. {\sl Triangles remarquables}

Dans la cat\'egorie $\dg(k)$ des motifs g\'eom\'etriques, nous utiliserons les triangles remarquables ([V-TCM]) suivants:

(1.3.5.1) $Gysin$.  Si Z est un sous-sch\'ema ferm\' e lisse, partout de codimension $c$, d'un sch\' ema lisse X,  
$$ \M(X-Z) \rightarrow \M(X)\rightarrow \M(Z)(c)[2c]\rightarrow\M(X-Z)[1].$$
(1.3.5.2) $Localisation\ \grave a\ support\ compact$. Si $Z$ est un sous-sch\'ema ferm\'e de $X$,
$$ \M^c(Z) \rightarrow \M^c(X) \rightarrow \M^c(X-Z)\rightarrow \M^c(Z)[1].$$
(1.3.5.3) $ Gysin\ g\acute en\acute eralis\acute e$. Si $X$ est un sch\'ema lisse \'equidimensionnel de dimension $n$ et $Z$ est un sous-sch\'ema ferm\'e de $X$,
$$\M(X-Z)\rightarrow \M(X)
\rightarrow \M^c(Z)^*(n)[2n] \rightarrow \M(X-Z)[1].$$
(1.3.5.4) $Eclatement$.  Si $Z$ est un sous-sch\'ema ferm\'e d'un sch\'ema $X$ et $p_Z : X_Z \rightarrow X$ est \'eclatement de $Z$ dans $X$,
$$\M(p_Z^{-1}(Z)) \rightarrow\M(X_Z)\oplus \M(Z)\rightarrow  \M(X) \rightarrow \M(p_Z^{-1}(Z))[1].$$
On remarquera que le triangle de Gysin g\'en\'eralis\'e est obtenu par dualit\'e \`a
partir du triangle de localisation \` a support compact.

\beginsection{1.4. Changements de base}

\th PROPOSITION 1.4.1.
\enonce Toute extension s\'eparable de corps $\sigma : k \hookrightarrow K$ induit un foncteur de changement de base
$$\matrix{
\sigma^s:& \sh_{Nis} (\smc(k)) & \rightarrow & \sh_{Nis} (\smc(K)) \cr
 & F & \mapsto & F_K\cr
}$$
envoyant, pour tout sch\'ema lisse $X$, le faisceau $\Ztr(X)$ sur le faisceau $\Ztr(X_K)$ avec $X_K =  X\times_{\sp(k)} \sp(K)$.

Le foncteur $\sigma^s$ est exact.
\endth

{\pc DEMONSTRATION}: le foncteur extension des scalaires $\alpha_K: \sm(k) \rightarrow \sm(K)$, d\'efini par $X\mapsto X\times_{\sp(k)} \sp(K)$ respecte les correspondances finies ([MVW] 1.12), est continu et permet de construire un foncteur changement de base. Par construction, on a pour tout faisceau avec transferts $F$ sur $\sm(k)$ et tout sch\'ema $Y$ lisse sur $K$ 
$$F_K(Y) = \limind_{(X,f)} F(X)$$
o\` u la limite est prise sur la cat\'egorie des couples $(X,f)$, o\`u $X$ est un sch\'ema lisse sur $k$ et $f$ une correspondance de $Y$ vers $X_K$. Comme l'extension est s\'eparable, un tel $Y$ est limite projective filtrante de vari\'et\'es lisses de type fini sur $k$, reli\'ees par des morphismes affines: le foncteur $F\mapsto F_K$ est exact. 
\hfill\cqfd

\smallskip
Puisque le foncteur $\sigma^s$  pr\'eserve l'invariance d'homotopie et le motif de Tate, qu'il est compatible au produit, il induit un foncteur de changement de base
$$\matrix{
\sigma_K :&\dm(k) & \rightarrow & \dm(K) \cr
 & \M & \mapsto & \M_K\cr
}$$
envoyant pour tout sch\'ema $X$ le motif $\M(X)$ sur le motif $\M(X_K)$.

\smallskip
{\pc REMARQUE} 1.4.2. Si l'extension  $\sigma: k \hookrightarrow K$ est finie, alors $\sp K$ d\'efinit un objet de $\dm(k)$ et le foncteur de changement de base admet un adjoint \`a gauche $\sigma_K^\ast$ induit par le foncteur de
$\sm(K)$ dans $\sm(k)$, qui \`a $Y$ associe $Y$; plus pr\'ecis\'ement, si $\M$ est un objet de $\dm(k)$, alors $\sigma_K^\ast\circ\sigma_K(\M)\simeq \M\otimes \M(\sp K)$, que l'on note plus simplement $\M\otimes K$.

\smallskip
En [LW09] nous avons d\'efini pour tout corps de caract\'eristique $0$ un ind-motif de De Rham, limite inductive de complexes motiviques. Nous avons besoin du r\'esultat suivant
\th LEMME 1.4.3.
\enonce Si $k \hookrightarrow K$ est une extension de corps de caract\'eristique $0$, le foncteur changement de base respecte les motifs de De Rham.
\endth
{\pc DEMONSTRATION}: le foncteur de changement de base respectant les faisceaux avec transferts, il suffit de v\'erifier que le foncteur de changement de base 
 $\widetilde{\sm k}_{Ab}\rightarrow  \widetilde{\sm K}_{Ab}$ envoie le faisceau des $k$-diff\'erentielles de K\"ahler sur le faisceau des $K$-diff\'erentielles de K\"ahler. Localement, c'est la formule de changement de base des diff\'erentielles ([EGA IV] 16.6.4).
\hfill\cqfd

\beginsection{1.5. La filtration par les poids de Bondarko}

Utilisant la construction de la cat\'egorie des motifs g\'eom\'etriques \` a partir de la cat\'egorie des complexes born\'es de sch\'emas projectifs de $\smc(k)$, Bondarko munit la \cat\ $\dg(k)$ d'une structure diff\'erentielle gradu\'ee [Bo09], la graduation sur les morphismes \'etant induite par celle du complexe de Suslin, plus pr\'ecis\'ement par le complexe cubique de Suslin ([Bo09] Ch.1).
Cela lui permet de munir $\dg(k)$ d'une structure  \` a poids ou pond\'erale ("weight structure"), \` a savoir 

\th THEOREME 1.5.1. ([Bo10] 1.1.1. et 6.5.3.)

\enonce Il existe deux sous-cat\'egories $\D^{w\geq 0}$ et $\D^{w\leq 0}$ de $\dg(k)$ telles que
{\parindent=1cm
\item{(i)} $\D^{w\geq 0}$ et $\D^{w\leq 0}$ sont additives et Karoubiennes;
\item{(ii)} {\it semi-invariance par translation:} $\quad\D^{w\geq 0}\subset\D^{w\geq 0}[1]$ et $\D^{w\leq 0}[1] \subset \D^{w\leq 0}$;
\item{(iii)}{\it orthogonalit\'e:} pour tout  objet $\M$ de $\D^{w\geq 0}$ et tout  objet ${\bf L}$ de $\D^{w\leq 0}[1]$, on a 
$$\ho_{\dg(k)} (\M, {\bf L}) = \lbrace0\rbrace;$$
\item{(iv)}{\it d\'ecomposition en poids:} pour tout motif g\'eom\'etrique $\M$, il existe un  triangle distingu\'e
$$ \M \rightarrow {\bf A}  \rightarrow {\bf B}\rightarrow \M[1]$$
avec $ {\bf A}$ objet de $\D^{w\leq 0}$ et ${\bf B}$ objet de $\D^{w\geq 0}$.
\par} 
\endth

{\pc REMARQUES}:

1.5.2. Les cat\'egories $\D^{w\geq 0}$ et $\D^{w\leq 0}$ sont construites respectivement \`a partir des classes de complexes born\'es de sch\'emas projectifs lisses de $\smc(k)$ qui sont concentr\'es en degr\'es positifs ou n\'egatifs. Plus pr\'ecis\'ement, si l'on note $\D^{w=0}=\D^{w\leq 0} \cap \D^{w\leq 0}$, la cat\'egorie des motifs de poids $0$, Bondarko l'identifie \`a la cat\'egorie de Grothendieck des motifs de Chow ([Bo10] 6.2) en utilisant les r\'esultats de Friedlander et Voevodsky [FV]. Ainsi les sommes finies de motifs de sch\'emas projectifs lisses et leurs facteurs directs, comme les $\Z(q)[2q]$, sont de poids $0$. Sur un corps qui v\'erifie la r\'esolution des singularit\'es, on montre en utilisant les triangles de Gysin que les sch\'emas lisses sont de poids positifs ([Bo09] 6.2.1).

\smallskip
1.5.3. La d\'ecomposition en poids (iv) n'est pas unique. N\'eanmoins Bondarko montre que l'on peut fixer a priori pour chaque motif g\'eom\'etrique $\M$ des d\'ecompositions
$$\M[i]\rightarrow \M^{w\leq i} \rightarrow \M^{w\geq i+1} \rightarrow \M[i+1]$$
et construire un complexe ([Bo10] 2.2) dit complexe des poids 
$$\cdots  \M^{(i-1)}\hfl{p_{i-1}}{} \M^{(i)}\hfl {p_i}{} \M^{(i+1)}\rightarrow\cdots$$avec ([Bo10] 1.5.6)
$$\matrix{
\M^{(i)} & \simeq & {\rm c\hat one} & (& \M^{w\leq i} [-1] & \rightarrow & \M^{w\leq i-1}&)&\cr
       &\simeq & {\rm c\hat one} & (& \M^{w\geq i+1} [-1] & \rightarrow & \M^{w\geq i}&)&\cr
}$$
et $\M^{(i)}$ est un motif de poids $0$.  
Un tel complexe n'est pas unique et sa construction n'est donc pas fonctorielle, mais, lorsqu'on lui applique un foncteur cohomologique, il fournit une suite spectrale qui est fonctorielle \`a partir de $E_2$ (cf ci-dessous 1.5.9.).

\smallskip
Il est possible de donner une description explicite d'un complexe des poids du motif ${\bf M}(X)$, lorsque $X$ est un sch\'ema  lisse et quasi-projectif sur $k$. Les notations sont celles de [D71], Ch.3.

\th{PROPOSITION} 1.5.4.
\enonce Soit $X$ un sch\'ema  lisse et quasi-projectif sur $k$, plong\'e  
dans un sch\'ema projectif lisse ${\bar X}$, tel que le sch\'ema compl\'ementaire $Y = {\bar X}-X$ soit un diviseur \`a croisements normaux et \`a composantes irr\'eductibles lisses $Y = \cup_{i=0}^{N-1} Y_i$. D\'esignons par $Y^j$ (resp. $\widetilde {Y^j}$) la r\'eunion (resp. somme disjointe) des intersections $j$ \`a $j$ des $Y_i$. On pose $\widetilde {Y^0}= Y^0 = {\bar X}$ et $\widetilde Y=\widetilde {Y^1}= \coprod_{0\leq i \leq N-1} Y_i$.

Le complexe $(\M(X)^{(j)})$ tel que, en degr\'e $j$,
$$\M(X)^{(j)} = \left\{\matrix{\M (\widetilde {Y^j})(j)[2j] & \hbox{ si }\   0\leq j \leq N\cr
0\hfill & \hbox{ sinon }\cr
}\right.$$
et dont les diff\'erentielles sont les sommes altern\'ees des morphismes induits par les inclusions de $\widetilde {Y^{j+1}}$dans $\widetilde {Y^j}$ est un complexe des poids de $\M(X)$.
\endth
\goodbreak
{\pc DEMONSTRATION}: on a, pour l'inclusion  $X\hookrightarrow {\bar X}$,  le triangle de Gysin g\'en\'eralis\'e, o\`u $n={\rm dim}(X)$
$$ \M(X) \rightarrow \M({\bar X}) \rightarrow \M(Y)^*(n)[2n] \rightarrow \M(X) [1].
\leqno (1.5.5)$$
Le motif $\M({\bar X})$ est de poids $0$ et si l'on prouve que le motif $\M(Y)^*(n)[2n]$
est de poids positif, alors le triangle (1.5.5) est une d\'ecomposition en poids de $\M(X)$ et des manipulations \'el\'ementaires dans les cat\'egories pond\'erales ([Bo10], lemme 1.5.4) prouvent que l'on peut choisir comme complexe des poids de $\M(X)$  un complexe des poids de 
$\M(Y)^*(n)[2n]$ d\'ecal\'e d'un degr\'e et augment\'e de $\M({\bar X})$ en degr\'e $0$. Plus pr\'ecis\'ement, on obtient
$$\eqalign{ \M(X)^{(j)} &= \M (Y)^*(n)[2n])^{(j-1)}\quad  {\rm si}\   1\leq j \cr
\M(X)^{(0)} &= \M(\bar X).}$$
 Ainsi il est \'equivalent de calculer les poids de $\M(X)$ ou de $\M(Y)^*(n)[2n]$.
Tout repose sur le lemme
\th{LEMME} 1.5.6.
\enonce Soit $Q = \cup_{i=0}^{N-1} Q_i$ un diviseur \`a croisements normaux, \`a $N$ composantes irr\'eductibles lisses, d'une vari\'et\'e lisse $P$ de dimension $n$. Un complexe des poids de 
$\M=\M(Q)^*(n)[2n]$ est 
$$ \M^{(i)}= \left\{\matrix{
\M(\widetilde{Q^{i+1}})(i+1)[2i+2] & {\rm si}\ 0\leq i\leq N-1\cr
0\hfill & {\rm sinon}\cr
}\right.$$
\endth
{\pc DEMONSTRATION}: on proc\`ede par r\'ecurrence sur le nombre $N$ de composantes de $Q$. Si $N=1$, le sch\'ema $Q$ est projectif lisse de dimension $n-1$ et le motif 
$\M(Q)^*(n)[2n]= \M(Q)(1)[2]$ est de poids $0$. 

Supposons le lemme d\'emontr\'e pour $N$ composantes et consid\'erons le sch\'ema 
$Q= \cup_{i=0}^N Q_i\hookrightarrow P$. Soit $Z = \cup_{i=0}^{N-1} Q_i$ le diviseur \`a croisements normaux \`a $N$ composantes. Le triangle de localisation \`a support compact (1.3.5.2) pour l'inclusion 
$Z\hookrightarrow Q$ est
$$\M(Z) \rightarrow \M(Q)\rightarrow \M^c(Q-Z)\rightarrow \M(Z)[1].\leqno(1.5.7)$$
Le sch\'ema $Q-Z= Q_N -(Q_{N0}\cup \dots \cup Q_{NN-1})$, avec $Q_{Ni}= Q_N\cap Q_i$, est lisse et son complexe des poids est connu par r\'ecurrence.
Par dualit\'e, on a 
$$\M^c(Q-Z)^* =\M(Q-Z)(-n+1)[-2n+2].$$
Le dual du triangle (1.5.7) devient apr\`es torsion par $\Z(n)[2n]$
$$\M(Q-Z)(1)[2] \rightarrow \M\rightarrow \M(Z)^*(n)[2n]\rightarrow \M(Q-Z)(1)[3].\leqno(1.5.8)$$
En appliquant \`a nouveau le lemme 1.5.4 de [Bo10], le triangle (1.5.8.) fournit en chaque poids un triangle distingu\'e qui est scind\'e, car chacune des extr\'emit\'es est somme de motifs purs de poids $0$ ([Bo10] Prop. 1.3.1.(7)). Par cons\'equent 
le complexe des poids de $\M$ est la somme des deux complexes des poids de $\M_1 = \M(Q-Z)(1)[2]$ et de $\M_2= \M(Z)^*(n)[2n]$. On a en degr\'e $i$ les intersections $(i+1)$ \`a $(i+1)$ des $Q_i$, tordues par $\Z(i+1)[2i+2]$; dans $\M_1$ viennent celles o\`u appara\^it $Q_N$, dans $M_2$ celles o\`u la composante $Q_N$ n'appara\^it pas.
\hfill\cqfd
\medskip
1.5.9. Pour tout foncteur cohomologique $H$ de $\dg(k)$ vers une cat\'egorie ab\'elienne $A$ et tout entier $i$, on pose
$$(W_i H)(\M) = {\rm Im} (H(w_{\leq i} \M ) \rightarrow H(\M))$$
o\`u l'on a not\'e $w_{\leq i} \M = \M^{w\leq i} [-i]$. 
On obtient ainsi une filtration croissante sur $H(\M)$, qui ne d\'epend pas du choix de $\M^{w\leq i}$, appel\'ee {\sl filtration par le poids}. D'autre part, on dispose de la suite spectrale associ\'ee au complexe des poids de $\M$
$$E_1^{pq}(\M) = H^q (\M^{(-p)}) \Rightarrow H^{p+q} (\M).$$
Le terme $E_1$ d\'epend du choix du complexe des poids, mais pas le terme $E_2$ et la suite est fonctorielle en $\M$ \`a partir de $E_2$ (cf [Bo10], Th 2.4.1).

Pour un motif $\M$, notons $\gr^W_iH(\M)=(W_i H)(\M)/(W_{i-1} H)(\M)$; on constate que 
$$\eqalign{(W_iH)(\M[-p])&={\rm Im} (H(w_{\leq i}(\M[-p]) ) \rightarrow H(\M[-p]))\cr
&={\rm Im} (H((w_{\leq i-p}\M)[-p] ) \rightarrow H(\M[-p]))\cr
&={\rm Im} (H^p(w_{\leq i-p}\M) \rightarrow H^p(\M))\cr
&=(W_{i-p}H^p)(\M)
}$$
d'o\`u le d\'ecalage habituel
$\gr^W_iH(\M[-p])=\gr^W_{i-p}H^p(\M)$ (cf [D71]). 

Une transformation naturelle entre deux foncteurs cohomologiques pr\'eserve les filtrations par le poids.

\beginsection{2. R\'ealisation de De Rham}

\beginsection{2.1. Construction}

En [LW09], nous avons muni les faisceaux $X \mapsto \Omega^n _{X/k} (X)$ de 
transferts et d\'efini un ind-complexe motivique  $\om$ qui 
repr\'esente la cohomologie de De Rham, dans le sens que pour tout sch\'ema 
$X$ lisse sur $k$, on a 
$${\Bbb H}^p_{Zar} (X,  \Omega^\cdot_{X/k} ) = 
\ho_{\DM(k)} (\M (X),\om [p])$$
o\`u $\Omega^\cdot _{X/k}$ est le complexe de De Rham de $X$ et 
l'on a plong\'e la cat\'egorie $\dm(k)$ dans la cat\'egorie des complexes motiviques non born\'es $\DM(k)$ de Cisinski et D\'eglise [CD09].
On g\'en\'eralise la cohomologie de de Rham \`a tous les motifs en posant

{\pc DEFINITION} 2.1.1. La {\sl r\'ealisation de De Rham} de tout complexe motivique $\M$ est le $k$-espace vectoriel gradu\'e
$$ {\bf H}^\cdot_{DR}(\M)= \oplus_{p\geq 0}{\bf H}^p_{DR}(\M)
= \oplus_{p\geq 0}{\bf H}_{DR}(\M [-p])$$
associ\'e au foncteur cohomologique 
${\bf H}_{DR} $ de $ \dm(k)$ dans la cat\'egorie des $k$-espaces vectoriels
$${\bf H}_{DR} (\M)  = \ho_{\DM (k)} (\M, \om).$$

\smallskip\goodbreak
Cette d\'efinition co\"\i ncide avec la d\'efinition de [LW09] o\` u nous tronquions le complexe de De Rham  en vertu du lemme suivant.

\th{LEMME} 2.1.2.
\enonce Pour tout complexe motivique $\M$ de $\dm(k)$, on a pour tout entier 
$i$, les isomorphismes 
$$\ho_{{\bf DM}(k)}(\M, \om[i]) \simeq \limind_n\ho_{\dm(k)} (\M,\tau_{\leq n}\om 
[i])$$
o\`u  $\tau_{\leq n}$ est le foncteur de filtration canonique \`a droite.
\endth

{\pc DEMONSTRATION}:  la propri\'et\'e ci-dessus \'etant stable par 
quasi-isomorphisme, d\'ecalage et c\^one, il suffit par le lemme 9.3 de [MVW]  de 
v\'erifier que pour toute famille de sch\'emas lisses 
$(X_{\alpha})_{\alpha \in A}$ on a 
$$ \prod_{\alpha \in A} \limind_n \ho_{D^-} (\Ztr (X_{\alpha}), \tau_{\leq n}
\om[i]) = \limind_n \prod_{\alpha \in A} \ho_{D^-} (\Ztr (X_{\alpha}), 
\tau_{\leq n}\om[i])$$
ce qui provient du fait que les complexes $\Ztr (X_{\alpha})$ sont concentr\'es 
en
 degr\'e z\'ero, les complexes  $\tau_{\leq n}\om[i]$ sont born\'es et la suite
  de droite est stationnaire pour $n>i$.
\hfill\cqfd
\goodbreak
\smallskip
Comme la cohomologie de de Rham des sch\'emas projectifs lisses est de dimension finie, le principe (0.2) implique que, restreinte aux motifs g\'eom\'etriques, la r\'ealisation de De Rham est un $k$-espace vectoriel de dimension finie.

\th{LEMME} 2.1.3.
\enonce Le morphisme de faisceaux de Nisnevich $d\!\log : {\cal O}^* \rightarrow \Omega^1$ commute aux transferts.
\endth

{\pc DEMONSTRATION}: en [LW09], nous avons d\'efini le transfert sur le faisceau des 
diff\'erentielles $\Omega^1$ en passant aux diff\'erentielles de Zariski $\Omega^{Zar}$, au sens de [K73].
Soit $Z$ une correspondance d'un sch\'ema lisse irr\'eductible $X$ vers un sch\'ema lisse
$Y$. Comme chez Suslin et Voedvodsky [SV96] on se ram\`ene au cas o\`u $Z$ est la normalis\'ee de $X$ dans une extension galoisienne finie du corps des fonctions $K(X)$ de $X$, de groupe de Galois $G= \Gal(K(Z)/K(X))$. On doit d\'emontrer la commutativit\'e du diagramme suivant:
$$\matrix{
{\cal O}^*_Z(Z) & \hfl{d\!\log}{} &  \Omega^1_Z(Z) \cr
\vfl{}{N} && \vfl{}{T_{Z/X}}\cr
{\cal O}^*_X(X)& \hfl{d\!\log}{}& \Omega^1_X(X)\cr
}$$
o\`u $N$ est la norme $N(f) = \prod_{\sigma \in G} \sigma (f)$ et $T_{Z/X}$ le transfert qui a \'et\'e d\'efini en (loc. cit) comme la compos\'ee
$$\Omega^1_Z(Z)\hfl{\alpha_Z}{}\Omega^{Zar}_Z(Z) \hfl{\sum_{\sigma \in G} \sigma^*}{}
\Omega^{Zar}_X(X) \hfl{\alpha_X^{-1}}{\simeq} \Omega^1_X(X).$$
Ici, pour tout sch\'ema $U$, on d\'esigne par $\Omega^{Zar}_U$ le faisceau bidual (au sens de ${\cal O}_U$-module dans la topologie de Zariski sur $U$) de $\Omega^1_U$ et
$\alpha_U: \Omega^1_U\rightarrow \Omega^{Zar}_U$ l'application canonique qui consiste \`a quotienter par la torsion. Le morphisme $\alpha_U$ est un isomorphisme quand le sch\'ema $U$ est lisse.

Si le sch\'ema $Z$ est lisse, le transfert co\" \i ncide avec la trace 
$Tr = \sum_{\sigma \in G} \sigma^*$ et le lemme est la traduction de la propri\'et\'e
$d\!\log \circ N = Tr \circ d\!\log$.

Si le sch\'ema $Z$ n'est que normal, il suffit de v\'erifier que cette propri\'et\'e n'est pas alt\'er\'ee par le passage au bidual. En remarquant que pour toute fonction inversible $f$ de ${\cal O}^*_Z(Z)$, la forme $\omega = d\!\log f$ v\'erifie la propri\'et\'e caract\'eristique des 
diff\'erentielles de Zariski sur un sch\'ema normal, \`a savoir 
$$P(Z): \forall x\in Z, \forall D \in Der_k({\cal O}_x, {\cal O}_x), \tilde D (\omega) \in {\cal O}_x$$ 
o\`u $Der_k({\cal O}_x, {\cal O}_x)$ est l'espace des $k$-d\'erivations de l'anneau local ${\cal O}_x$ dans 
lui-m\^eme et $\tilde D : \Omega^1_{{\cal O}_x/k}\rightarrow {\cal O}_x$ est le morphisme ${\cal O}_x$-
lin\'eaire canoniquement associ\'e \`a la diff\'erentielle $D$.

Or on a, pour toute fonction $f\in {\cal O}^*_Z(Z)$ et tout point $x$ de $Z$,
$$\tilde D (d\!\log f) = \tilde D ({df_x \over f_x}) = {Df_x\over f_x} \in {\cal O}_x.$$
\hfill\cqfd

\goodbreak
\smallskip
Ce morphisme $d\!\log$ induit le morphisme de  $\dm(k)$

$$\matrix{
\Z(1)&: &0& \hfl{}{} & 0      &\hfl{ }{}&{\cal O}^*   & \hfl{}{}& 0 \cr
            & & &        &        &         & \vfl{d\!\log}{}&         &  \cr
\tau_{\leq 2}\om&:&0& \hfl{}{} &{\cal O}&\hfl{d}{}& 
\om^1&\hfl{d}{}&{\rm Ker}(d)\cr
}\leqno (2.1.4)$$
qui produit un g\'en\'erateur de    ${\bf H}^0_{DR} (\Z (1) )$.
On v\'erifie que l'image de $d\log$ se trouve dans $\ke d$, d'o\`u un morphisme ${\bf Z}(1)\rightarrow \tau_{\leq 1}\om$.

 Un  g\'en\'erateur de ${\bf H}^0_{DR} (\Z (q) )$, lorsque $q$ est un entier positif, est induit par le morphisme
 issu du produit sur $\om$ [LW09]
 $$\Z(q) = \Z(1)^{\otimes q}\buildrel {d\!\log^{\otimes q}} 
 \over{\longrightarrow}(\om)^{\otimes q}\rightarrow \om.\leqno{(2.1.5)}$$
  
Ces morphismes permettent de d\'efinir des classes de Chern

$$c^{p,q}_{DR} : H^{p,q}(\M) = \ho_{\dm (k)} (\M, \Z(q)[p])
\rightarrow {\bf H}^p_{DR}(\M)$$
 pour tout complexe motivique $\M$ de $\dm$.

Le produit de [LW09] fournit aussi pour toute paire de complexes motiviques  $\M_1$ et
$\M_2$ de $\dm(k)$ un accouplement
$$ {\bf H}^{p_1}_{DR} (\M_1) \otimes {\bf H}^{p_2}_{DR} (\M_2) \rightarrow 
{\bf H}^{p_1+p_2}_{DR}(\M_1 \otimes\M_2).$$

\smallskip\goodbreak

\beginsection{2.2. Filtration par les poids}

La filtration par les poids est induite sur le foncteur cohomologique par la structure de poids d\'efinie par Bondarko [Bo10]
$$W_i {\bf H}_{DR}(\M) ={\rm Im} ({\bf H}_{DR}(w_{\leq i} \M ) \rightarrow {\bf H}_{DR} (\M)),$$
qui induit la filtration sur ${\bf H}^p_{DR}(\M)$
$$W_i {\bf H}^p_{DR}(\M) = W_{i+p}{\bf H}_{DR}(\M[-p]).\leqno(2.2.1)$$

\th{PROPOSITION} 2.2.2.
\enonce Si $X$ est un sch\'ema lisse quasi-projectif sur $k$, la filtration $W$ sur $ {\bf H}_{DR}^{\,\cdot}(\M(X))$ co\"\i ncide avec la filtration par le poids classique de la cohomologie de $X$.
\endth
\smallskip\goodbreak
{\pc DEMONSTRATION}: nous choisissons comme r\'ef\'erence les travaux de Deligne en th\'eorie de Hodge, plus pr\'ecis\'ement ([D71] Ch.3), repris sur un corps $k$ par Jannsen ([J90]). Les calculs, \`a base de p\^oles logarithmiques, y sont faits en cohomologie analytique mais s'adaptent en cohomologie de De Rham. 

Reprenons les notations de la proposition 1.5.4: 
le sch\'ema $X$ \'etant lisse et quasi-projectif sur $k$ comme dans (loc. cit)(3.2.1), on le plonge 
dans un sch\'ema projectif lisse ${\bar X}$, tel que le sch\'ema compl\'ementaire $Y = {\bar X}-X$ soit un diviseur \`a croisements normaux et \`a composantes irr\'eductibles lisses $Y = \cup_{i=0}^{N-1} Y_i$. Comme en (loc. cit.), on d\'esigne par $Y^j$ (resp. $\widetilde {Y^j}$) la r\'eunion (resp. somme disjointe) des intersections $j$ \`a $j$ des $Y_i$. On pose $\widetilde {Y^0}= Y^0 = {\bar X}$ et $\widetilde Y=\widetilde {Y^1}= \coprod_{0\leq i \leq N-1} Y_i$. On a vu que $\M(X)^{(0)}=\M(\bar X)$ et $\M(X)^{(k)}=\M(\widetilde{Y^k})(k)[2k]$ pour $1\leq k\leq N-1$.

La suite spectrale de Bondarko ([Bo10] 2.4.1.)
$$E_1^{pq}(\M(X)) = {\bf H}_{DR}^q (\M(X)^{(-p)}) \Rightarrow {\bf H}_{DR}^{p+q}(\M(X))$$
co\"\i ncide alors avec la suite spectrale de Jannsen ([J90], I3.) ou de Deligne ([D71] 3.2.7) lorsque $k=\bf C$
$${}_W E_1^{p,q}(X) = {\Bbb H}^{q+2p} (\widetilde {Y^{-p}}, \Omega^\cdot_{\widetilde{Y^{-p}}})\otimes_k {\bf H}^0_{DR}({\bf Z}(-p)).$$
On en d\'eduit qu'elle d\'eg\'en\`ere en $E_2$ (cf. corollaire 3.2.13, [D71], lorsque $k=\bf C$) ainsi que l'isomorphisme 
$$E_2^{p,q}(\M(X))=\gr^W_q{\bf H}_{DR}^{p+q}(\M(X)),$$
qui tient compte du d\'ecalage.\hfill\cqfd

\beginsection{2.3. Filtration de Hodge}

Dans la cat\'egorie d\'eriv\'ee $D(\sh_{Nis}(\smc(k)))$ des complexes de faisceaux avec transferts, la filtration b\^ete  ou troncature \`a gauche
$$(\sigma_{\geq p}\om)^q (X) =\left\{ \matrix{
\Omega^q_{X/k}(X) &{\rm  si}\ &p \leq q \cr
0\hfill &{\rm sinon}\  &\cr
}\right.$$
fournit pour tout complexe born\'e $K$ de $D^-_{Nis}(\smc(k))$ une suite spectrale
$$ E^{p,q}_1 (K) = \ho_{D(\sh_{Nis}(\smc(k)))}(K,\om^p[q]) \Rightarrow  \ho_{D(\sh_{Nis}(\smc(k)))}(K,\om[p+q]).\leqno(2.3.1)$$
Cette suite spectrale est non seulement fonctorielle en $K$ mais, de plus, comme pour la suite spectrale des poids de Bondarko, elle est fonctorielle sur $\dm(k)$ \` a partir du terme $E_2$: 
\th{PROPOSITION} 2.3.2.
\enonce Si $f:K\rightarrow K'$ est une ${\bf A}^1$-\'equivalence de complexes de faisceaux avec transferts dans $D^-_{Nis}(\smc(k))$, alors on a, pour tout $r\geq 2$ et tout couple d'entiers $(p,q)$, des isomorphismes
$$ E^{p,q}_r (K') \simeq E^{p,q}_r (K).$$
\endth
{\pc DEMONSTRATION}: il suffit de montrer que, pour tout couple d'entiers $(p,q)$, on a un isomorphisme 
$$E^{p,q}_2 (K') \simeq E^{p,q}_2 (K)$$
compatible 
 aux diff\'erentielles. En adaptant la d\'emonstration du lemme d'homotopie aux bicomplexes, on se ram\`ene \`a prouver que pour tout sch\'ema $X$ lisse sur $k$, la projection 
$\pi_X: X \times_k {\bf A}^1 \rightarrow X$ induit pour tout couple d'entiers $(p,q)$ un isomorphisme 
$$ E^{p,q}_2 ({\bf Z}_{tr}(X))\simeq  E^{p,q}_2 ({\bf Z}_{tr}(X\times_k {\bf A}^1)).\leqno (2.3.3)$$
Sur le site $X_{Nis}$, consid\'erons les complexes de faisceaux $\Omega_X^\cdot$ et 
$\pi_{X*}\Omega_{X\times_k{\bf A}^1}^\cdot$. Les membres de l'isomorphisme (2.3.3) sont les termes 
$E_2$ des suites spectrales associ\'ees \`a la filtration b\^ete des complexes respectifs. D'apr\`es [Ver96](4.4.3), ces suites co\"\i ncident \`a partir de $E_2$ lorsque les deux complexes sont homotopes.
Or la d\'ecomposition des formes diff\'erentielles et la formule de projection impliquent l'isomorphisme
$$\pi_{X*} \Omega^\cdot_{X\times {\bf A}^1} \simeq  \Omega^\cdot_X\otimes_{{\cal O}_X} \pi_{X*}\pi_{{\bf A}^1}^{-1}
\Omega^\cdot_{{\bf A}^1}$$
qui permet de remonter l'homotopie 
induite par l'int\'egration 
$$\matrix{
s:& \Omega^1_{{\bf A}^1}& \rightarrow & {\cal O}_{{\bf A}^1} \cr
& f(t) dt & \mapsto & \int_0^t f(u)du \cr
}$$
en une homotopie entre l'application compos\'ee $\Omega_{X\times_k{\bf A}^1}^\cdot\rightarrow \Omega_X^\cdot\rightarrow\Omega_{X\times_k{\bf A}^1}^\cdot$ et l'identit\'e.
\hfill\cqfd

\medskip
On a ainsi une suite spectrale $E^{p,q}_r (\M)$, $r\geq 2$ d\'efinie pour tout complexe motivique $\M$. Si de plus le motif $\M$ est g\'eom\'etrique, le principe (0.2) implique que les espaces vectoriels sont de dimension finie sur $k$ et que la suite spectrale (2.3.1) d'aboutissement ${\bf H}_{DR}(\M)$ d\'eg\'en\`ere et d\'efinit sur ${\bf H}_{DR}(\M)$ la filtration r\'eguli\`ere
$$F_{\rm naive}^q  {\bf H}^p_{DR}(\M) ={\rm Im}\left( \ho_{D(\sh_{Nis}(\smc(k)))}(\M,\sigma_ {\geq q}
\om[p])\rightarrow \ho_{\DM}(\M,\om[p])\right).
\leqno (2.3.4)$$
On a par construction 

\th{PROPOSITION} 2.3.5. 
\enonce Si $\M$ est le motif associ\'e \`a un sch\'ema $X$ lisse et projectif, la filtration $F_{\rm naive}$ sur ${\bf H}_{DR}(\M)$ est la filtration de Hodge de la cohomologie du sch\'ema $X$. 
\endth
\smallskip
Si le sch\'ema lisse $X$ n'est pas projectif, c'est traditionnellement en rempla\c cant le complexe de de Rham par un complexe de de Rham \`a p\^oles logarithmiques, qui permet aussi de retrouver les poids, que l'on d\'ecrit la filtration de Hodge. De fa\c con similaire, nous r\'esolvons par les poids en rempla\c cant tout motif g\'eom\'etrique par un complexe des poids pour obtenir la filtration de Hodge.

\th{PROPOSITION} 2.3.6.
\enonce Soient $\M$ un motif g\'eom\'etrique effectif, $\M^{(\cdot)}$ un complexe des poids associ\'e \`a $\M$. La cohomologie du complexe total ${\bf K}_{DR}$ du bicomplexe $\r^\cdot\ho (\M^{(\cdot)},{\bf\Omega})$ calcule la r\'ealisation de De Rham de $\M$.
\endth
{\pc DEMONSTRATION}

Pour un motif g\'eom\'etrique effectif, on peut construire un complexe des poids born\'e et quitte \`a translater, on se ram\`ene au cas o\`u $\M$ est de poids positif ou nul. Tout complexe des poids associ\'e \`a $\M$ s'\'ecrit $\M^{(0)}\rightarrow\cdots\rightarrow \M^{(M)}\rightarrow0$ avec $\M^{(0)}\not=0$.

La fl\`eche $\M \rightarrow \M^{w\leq 0}$ et l'\'egalit\'e $\M^{(0)}=\M^{w\leq 0}$ induisent une fl\`eche entre les complexes 
$${\bf K}_{DR}\rightarrow \r^\cdot\ho (\M,{\bf\Omega}).\leqno{(2.3.7)}$$

Le complexe  
$\r^\cdot\ho (\M,{\bf\Omega})$ est muni de la filtration par les poids \`a la Bondarko
$$W_k \r^\cdot\ho (\M,{\bf\Omega})= {\rm Im}(\r^\cdot\ho (w_{\leq k}\M,{\bf\Omega})\rightarrow \r^\cdot\ho (\M,{\bf\Omega}))$$
 qui fournit une suite spectrale ${}^WE_1^{p,q}(\r^\cdot\ho (\M,{\bf\Omega}))={\bf H}_{DR}^q(\M^{(-p)})$ convergeant fortement vers ${\bf H}_{DR}(\M)$.
D'autre part, la premi\`ere filtration du bicomplexe  $K_{DR}^{p,q}=\r^q\ho (\M^{(-p)},{\bf\Omega})$ fournit la suite spectrale $${}^WE_1^{p,q}({\bf K}_{DR})={\bf H}_{DR}^q(M^{(-p)}).$$ 
Puisque l'on peut choisir le complexe $\M^{(0)}\rightarrow\cdots\rightarrow \M^{(k)}\rightarrow0$ comme complexe des poids de $w_{\leq k}\M$, la fl\`eche $(2.3.7)$ est compatible aux filtrations et induit un isomorphisme entre les suites spectrales associ\'ees.

Le th\'eor\`eme 3.2, chapitre XV de [CE], permet de conclure.\hfill\cqfd

\smallskip
{\pc DEFINITION} 2.3.8. Soit $\M$ un motif g\'eom\'etrique effectif, $\M^{(\cdot)}$ un complexe des poids associ\'e. La {\sl filtration de Hodge} de ${\bf H}_{DR}(\M)$ est induite par la filtration du bicomplexe
$$F^j(\r^q\ho (\M^{(-p)},{\bf\Omega}))=\im(\r^q\ho (\M^{(-p)},{}_{\sigma\geq j}{\bf\Omega})\rightarrow \r^q\ho (\M^{(-p)},{\bf\Omega})),$$
et
$$F^j\r^\cdot\ho(\M,\om)=\im\left(F^j\r^\cdot\ho(\M^{(\cdot)},\om)\rightarrow\r^\cdot\ho(\M,\om)\right).$$

\smallskip
{\pc REMARQUE}: nous verrons en (4.2.4) que ces d\'efinitions ne d\'ependent pas des choix, notamment de celui du complexe des poids.

\beginsection{2.4. Motif de De Rham analytique}

En munissant la cat\'egorie $\smc(\C)$ des sch\'emas alg\'ebriques complexes lisses et des correspondaces finies de la topologie analytique, on fait les m\^emes constructions qu'en topologie  \'etale pour d\'efinir une cat\'egorie triangul\'ee $\dman$. Rappelons que le th\'eor\`eme d'invariance d'homotopie (1.3.3.1) n'est plus valable dans ce cadre. Le changement de topologie d\'efinit un foncteur 
$$\matrix{
t_{an}:& \dm(\C) &\rightarrow &\dman \cr
& \M & \mapsto & \M_{an}
}$$
D'apr\`es [GAGA], ce foncteur  envoie le faisceau structural $X \mapsto 
\Gamma (X,{\cal O}_X)$ sur le faisceau structural des fonctions holomorphes et 
les diff\'erentielles de K\"ahler sur les diff\'erentielles holomorphes, le ind-motif $\om$ vers le ind-motif not\'e $\om_{an}$ dont le $n^{\rm i\grave eme}$-terme est le faisceau $X\mapsto \Omega_X^{n,h}(X)$ des $n$-diff\'erentielles 
holomorphes. 
\th{PROPOSITION} 2.4.1
\enonce Le foncteur cohomologique $\M \mapsto \ho_{\DMan}(\M_{an} , \om_{an})$  de la cat\'egorie des motifs g\'eom\'etriques vers la cat\'egorie des $\bf C$-espaces vectoriels co\"\i ncide avec la r\'ealisation de De Rham.
\endth
{\pc DEMONSTRATION}: Le foncteur $t_{an}$ induit pour tout motif $\M$ g\'eom\'etrique de $\dg(\C)$ un morphisme 
$$ \ho_{\DM(\C)}(\M , \om)\rightarrow \ho_{\DMan}(\M_{an} , \om_{an}).$$
Si $\M=\M(X)$ est le motif d'un sch\'ema lisse projectif, ce morphisme est un isomorphisme d'apr\`es [GAGA]. On conclut en rappelant que de tels $\M(X)$ engendrent $\dg(\C)$.
\hfill\cqfd

\goodbreak
\beginsection{3. R\'ealisation de Betti}

Cette partie est largement inspir\'ee des travaux de Suslin et Voevodsky [SV96] et a \'et\'e r\'esum\'ee en [L08].

\beginsection{3.1. R\'ealisation topologique des motifs sur $\C$.}

Nous travaillons sur le site $\sm (\C)$ des sch\'emas lisses sur le corps des nombres 
complexes, muni de la topologie de Nisnevich, et sur le site $CW$ des espaces 
topologiques r\'eels admettant une triangulation, muni de la topologie des 
hom\'eomorphismes locaux.
 
Le foncteur $\theta:\sm (\C)\rightarrow CW$ qui \`a un sch\'ema $X$ associe la vari\'et\'e $X(\C)$ des points complexes est continu. Il induit un foncteur $\theta^s : \widetilde{\sm \C } \rightarrow \widetilde{CW}$ prolongeant $\theta$.

\th{PROPOSITION} 3.1.1.
\enonce Si $X$ est un sch\'ema projectif lisse, l'image 
$\theta^s (\Ztr (X))$ est le faisceau 
$$U \mapsto {\rm Hom} (U, \coprod_{d\geq 0} S^d X(\C))^+ $$
o\`u le sch\'ema $S^dX$ est la puissance sym\'etrique du sch\'ema $X$ et pour tout mono\"\i de $M$,  on d\'esigne par $M^+$ le groupe groupe de Grothendieck associ\'e.
\endth
{\pc DEMONSTRATION}: le r\'esultat repose essentiellement sur le th\'eor\`eme de Suslin et Voevodsky.

\th THEOREME 3.1.2. (Suslin et Voedvosky [SV96 Theorem 6.8])

\enonce Si $X$ est un sch\'ema quasi-projectif sur un corps $k$, on a, pour tout
 sch\'ema (normal et connexe) $S$, un isomorphisme de groupes
$$ {\bf Z}_{tr} (X)(S) = \ho_{\sch k}(S,\coprod_{d\geq 0} S^d(X))^+\ ,$$
o\`u  la cat\'egorie $\sch k$ est celle des sch\'emas de type fini sur $k$.
\endth

Le foncteur $\theta^s$ respectant les foncteurs repr\'esentables, il envoie, pour tout entier positif $d$ et tout sch\'ema lisse $X$, le faisceau repr\'esent\'e par $X^d$ sur le faisceau repr\'esent\'e par $X^d(\C)$. Par ailleurs, le foncteur  $\theta^s$ commute aux colimites finies (comme le quotient par les groupes de permutation) mais \'egalement aux colimites quelconques ([SGA4] III 1.3): il envoie le faisceau repr\'esent\'e par le ind-sch\'ema lisse
$\coprod_{d\geq 0} S^dX$ sur le faisceau repr\'esent\'e par l'espace
$\coprod_{d\geq 0} S^d(X(\C))$. Comme le foncteur est compatible aux structures alg\'ebriques, la proposition (3.1.1.) s'en d\'eduit.
\hfill\cqfd

\bigskip
Composant le foncteur $\theta^s$ avec le foncteur exact oubli de transferts, nous obtenons un foncteur exact \`a droite qui se factorise dans les faisceaux ab\'eliens en un foncteur mono\" \i dal
$$\Phi : \sh_{Nis} (\smc \C) \rightarrow \widetilde{CW}_{Ab}$$
de la cat\'egorie des faisceaux de Nisnevich avec transferts vers le topos ab\'elien 
$\widetilde{CW}_{Ab}$.
En composant de plus avec le foncteur complexe de Suslin $C_*$ qui est exact, 
nous obtenons un foncteur, toujours exact \`a droite,
$$\Psi : \sh_{Nis} (\smc \C) \rightarrow C^-(\widetilde{CW}_{Ab})$$
o\`u $C^-(\widetilde{CW}_{Ab})$ est la cat\'egorie ab\'elienne des complexes 
born\'es sup\'erieurement de faisceaux ab\'eliens. C'est ce foncteur que nous d\'erivons en utilisant la classe $\Sigma$ des sommes de faisceaux $\oplus _{\alpha} \Ztr (V_{\alpha})$ o\`u 
les $V_{\alpha}$ parcourent les sch\'emas lisses quasi-projectifs tels que la vari\'et\'e 
des points complexes $V_{\alpha}(\C)$ soit \`a composantes connexes contractiles. Cette classe 
$\Sigma$ d'objets est bien s\^ur stable par somme directe et  permet via la proposition 1.3.1.1  de construire des r\'esolutions \`a gauche de tout faisceau avec transferts. Pour que le foncteur $\Psi$ admette un foncteur d\'eriv\'e \`a gauche, il suffit par [Ver77] de v\'erifier
que les objets de $\Sigma$ sont acycliques \`a gauche en montrant

\th LEMME 3.1.3.
\enonce Soit $N$ un complexe acyclique et born\'e sup\'erieurement d'objets 
$\oplus _{\alpha} \Ztr (V_{\alpha})$ de $\Sigma$. Alors $\Psi(N)$ est acyclique.
\endth

{\pc DEMONSTRATION }: pour tout entier $i$, nous devons montrer que le faisceau de Nisnevich 
associ\'e au pr\'efaisceau $H_i(\Psi (N))$ est nul. Cela revient \`a montrer que pour tout 
pr\'efaisceau $F$ tel que le faisceau associ\'e $F_{Nis}$ soit nul, les d\'eriv\'es 
$L^n\Phi(F)$ sont nuls. Comme en [MVW] 8.15 on se ram\`ene au cas d'un complexe 
de recouvrement de Cech et nous devons montrer le lemme suivant:
\smallskip 
\th LEMME 3.1.4.
\enonce Si ${\cal U}=\{U_i\rightarrow X\}$ est un recouvrement de Nisnevich d'un 
sch\'ema $X$ lisse quasi-projectif sur $\C$ par des ouverts dont les points complexes sont 
\`a composantes contractiles, alors l'image par le foncteur 
$\Phi$ du complexe de Cech $\check N ({\cal U}/X)$ est acyclique.
\endth
Ce lemme repose sur les travaux de Suslin et Voevodsky en [SV96] section 10 qui ont
r\'einterpr\'et\'e le th\'eor\`eme de Dold-Thom en 

\th PROPOSITION 3.1.5.
\enonce Si $X$ est un sch\'ema lisse quasi-projectif sur $\C$, l'image 
$\Phi (\Ztr (X))$ est quasi-isomorphe au complexe des cha\^\i nes 
singuli\`eres de la vari\'et\'e topologique $X(\C)$.
\endth 
{\pc DEMONSTRATION}: par dualit\'e il est \'equivalent de montrer que le 
bicomplexe de cocha\^\i nes 
$$Sing^. (X(\C)) \rightarrow Sing^.(\check N{\cal U} (\C))$$
est acyclique. D'apr\`es la proposition 3.1.5 le  bicomplexe 
$ Sing^.(\check N{\cal U} (\C))$ calcule la cohomologie de Cech $\check H^*(X(\C))$
de $X(\C)$ qui est \'egale \`a la cohomologie singuli\`ere de $X(\C))$ par un th\'eor\`eme 
de Cartan ([Go58], 5.9.2).
\hfill\cqfd

Nous en d\'eduisons un foncteur 
$L\Psi : D^- (\sh_{Nis} (\smc \C)) \rightarrow D(Ab)$
qui par le lemme d'homotopie (1.3.3.2) se factorise en un foncteur dit de r\'ealisation topologique 
$t_{\C} : \dme (\C) \rightarrow D ({\cal A} b)$.

\smallskip
Par le th\'eor\`eme d'Eilenberg-Zilber, le foncteur $\Psi$ est compatible au produit et son 
d\'eriv\'e $t_{\C}$ commute au produit  de $\dme (\C)$ \'egalement d\'efini \`a partir des 
r\'esolutions par les $\Ztr (X)$.

\medskip
L'image de  $\Z (1) = C_*\Ztr ({\bf G}_m^{\wedge 1})[-1]$ est le 
complexe calculant l'homologie singuli\`ere r\'eduite de $\C^*$, d\'ecal\'e de 
$-1$, d'o\`u $t_{\C} ( \Z (1)) \simeq 2i\pi\Z $, complexe concentr\'e en degr\'e $0$, par le th\'eor\`eme des r\'esidus. 
Comme le foncteur $t_{\C}$ est tensoriel et envoie le motif de Tate $\Z (1)$ sur un objet  inversible, il s'\'etend en un foncteur de $\dm(\C)$ et 
les r\'esultats de ce paragraphe se r\'esument en 
\th THEOREME 3.1.6.
\enonce Il existe un foncteur tensoriel de r\'ealisation topologique 
$$\matrix{
t_{\C} :& \dm (\C) & \rightarrow & D ({\cal A} b) \cr
& \M & \mapsto & \M(\C),\cr
}$$
qui pour le motif 
$\M =\M(X)$   associ\'e \`a un sch\'ema $X$ 
quasi-projectif lisse sur $\C$ 
permet de repr\'esenter la cohomologie singuli\`ere 
de $X(\C)$, c'est-\`a-dire
$$ H^p(X(\C), \Z ) =\ho _{D^- ({\cal A} b)} (t_{\C} (\M(X)), \Z [p]).$$
\endth

La 
compatibilit\'e du foncteur de r\'ealisation au produit impose
$$t_{\C} ( \Z (q)) = (2i\pi)^q \Z.$$

\bigskip
Comme, pour les sch\'emas alg\'ebriques sur $\C$, la topologie analytique est interm\'ediaire entre la topologie de Nisnevich et la topologie usuelle, le foncteur $t_{\C}$ se factorise \`a travers $t_{an}: \dm(\C) \rightarrow \dman$. Par le lemme de Yoneda, on a pour tout sch\'ema lisse $X$ et tout faisceau analytique avec transferts $F$, un isomorphisme 
$$\ho_{\dman} (\M(X)_{an}, F) = H_{an} (X(\C), F)$$
en particulier pour le faisceau constant $\Z$
$$\ho_{\dman} (\M(X)_{an}, \Z) = H_{sing} (X(\C), \Z).$$
On en d\'eduit, puisque les motifs de sch\'emas lisses engendrent $\dg(\C)$

\th PROPOSITION 3.1.7.
\enonce On a $t_{an} ( \Z (q)) = (2i\pi)^q \Z$ et la restriction de $t_{an}$ \`a $\dg(\C)$ induit des isomorphismes
$$\ho_{D^- ({\cal A} b)} (t_{\C}(\M), (2i\pi)^q\Z) \simeq \ho_{\dman} (\M(X)_{an}, \Z(q)_{an})$$
pour tout motif g\'eom\'etrique $\M$ sur $\C$.
\endth

\beginsection{3.2. R\'ealisation topologique des motifs sur $\R$.}

Si le sch\'ema $X$ est d\'efini sur $\R$, la vari\'et\'e analytique $X(\C) = (X 
\times_{\sp \R} \sp \C)(\C)$ est munie d'une action continue de la conjugaison 
complexe $F_{\infty}$. En suivant cette action dans la construction 
pr\'ec\'edente, 
on montre que le foncteur de r\'ealisation topologique se factorise en un 
diagramme
$$\matrix{
\dm (\R) & \hfl{t_{\C, F_{\infty}}}{} &  D ({\cal A} b^{\sigma_2}) \cr
\vfl{}{\otimes_{\R} \C} && \vfl{}{}\cr
\dm (\C)& \hfl{t_{\C}}{}& {D ({\cal A} b)}\cr
}$$
o\`u ${\cal A} b^{\sigma_2}$ est la cat\'egorie ab\'elienne des groupes 
ab\'eliens munis d'une involution et la fl\`eche de droite est induite par 
l'oubli de l'involution.

Le motif de Tate $\Z(1)$ est r\'eel et, sur sa r\'ealisation $ t_{\C} ( \Z (1))$, 
l'involution est induite par le changement d'orientation de ${\bf S}^1$ dans 
$\C^*$ et agit par multiplication par $-1$. 

\beginsection{3.3. R\'ealisation de Betti des motifs }

Pour tout plongement $\sigma : k\rightarrow \C$, l'extension 
des scalaires ${\bf\sigma_{\C}} : \dm (k) \rightarrow  \dm (\C)$ construite  en 
1.4, compos\'ee avec la r\'ealisation topologique, d\'efinit un foncteur
$$\matrix{
t_{\sigma} &:& \dm(k) &\rightarrow & D({\cal A}b) \cr
&&\M &\mapsto & \M_{\sigma} (\C) : = t_{\C} \circ {\bf\sigma_{\C}}(\M ). \cr
}$$
Si le plongement $\sigma$ est r\'eel, le foncteur $t_{\sigma}$ se factorise dans la 
cat\'egorie des groupes ab\'eliens munis d'une involution (3.2). 

On pose alors 
${\bf H}_\sigma(\M,q)=\ho_{D({\cal A}b)} 
(\M_{\sigma}(\C), (2i\pi)^q \Z)$ et l'on d\'efinit
\th DEFINITION 3.3.1.
\enonce Pour tout complexe motivique $\M$ de $\dm(k)$, tout plongement 
$\sigma:k\rightarrow \C$ et tout entier $q\geq 0$, la r\'ealisation enti\`ere (resp. 
r\'ealisation) de 
Betti ${\bf H}_{\sigma}^\cdot ( \ , \Z(q))$ (resp. ${\bf H}_{\sigma}^\cdot( \ ,q)$) est  le groupe 
ab\'elien 
(resp. $\Q$-espace vectoriel) gradu\'e sur $\Z$ 
$${\bf H}^\cdot_\sigma(\M,\Z(q))=\oplus_{p\in\bf Z}{\bf H}^p_{\sigma} (\M,\Z(q)),$$  et
$${\bf H}^p_{\sigma} (\M,q) =\Q \otimes_{\Z} {\bf H}^p_{\sigma} (\M,\Z(q)) .$$
Si le plongement $\sigma$ est r\'eel, le complexe $\M_{\sigma}(\C)$ est muni d'une 
involution induite par la conjugaison complexe et les r\'ealisations ${\bf H}^\cdot_{\sigma} ( \ , \Z(q))$ et ${\bf H}^\cdot_{\sigma}( \ ,q)$ h\'eritent de cette structure.
\endth

Le foncteur $t_{\sigma}$ induit directement des classes de Chern, pour des 
entiers $p$ et $q$ et tout complexe motivique $\M$ de $\dm(k)$

$$c^{p,q}_{\sigma} :  H^{p,q}(\M) = \ho_{\dm (k)} (\M, \Z(q)[p])
\rightarrow\ho_{D({\cal A}b)} (\M_{\sigma}(\C), (2i\pi)^q \Z [p])= 
{\bf H}^p_{\sigma}(\M, \Z(q)).$$

Si le complexe motivique $\M= \M(X)$ est le motif d'un sch\'ema quasi-projectif 
lisse sur $k$, les groupes de 
r\'ealisation de Betti co\"\i ncident d'apr\`es (3.1.5) avec les groupes de cohomologie 
singuli\`ere 
$${\bf H}^p_{\sigma} (\M(X),\Z(q)) = H^p(X(\C), (2i\pi)^q\Z).$$

3.3.2. {\sl Filtration par le poids} 

Par ailleurs, la r\'ealisation de Betti ${\bf H}_{\sigma}(\ ,{\bf Z})$, pour $q=0$, \'etant un foncteur cohomologique, elle h\'erite de la filtration par le poids d\'ecrite en (1.5.9), \`a savoir, pour $i\in\bf Z$,
$$W_i{\bf H}_{\sigma}({\bf M} ,{\bf Z})=\im({\bf H}_{\sigma}(w_{\leq i}{\bf M} ,{\bf Z})\rightarrow {\bf H}_{\sigma}({\bf M} ,{\bf Z})).$$
Le calcul fait en (2.2) dans le cadre de la r\'ealisation de De Rham montre que cette filtration co\" \i ncide avec celle de Deligne apr\`es avoir tensoris\'e par $\bf C$. En particulier, pour un sch\'ema lisse $X$ que l'on plonge dans  ${\bar X}$,  la suite spectrale de Bondarko associ\'e \`a ${\bf H}_{\sigma}$ est isomorphe (\`a renum\'erotation pr\`es) \`a la suite spectrale de Leray de la cohomologie singuli\`ere pour l'inclusion $X(\C) \hookrightarrow {\bar X}(\C)$. 
\smallskip
{\pc DEFINITION} 3.3.3. Pour tout $q\in\bf Z$ et tout motif g\'eom\'etrique $\bf M$, on appelle {\sl filtration par le poids} de la r\'ealisation de Betti ${\bf H}_{\sigma}( {\bf M},{\bf Z}(q))$, la filtration croissante index\'ee par $\bf Z$ d\'eduite par la filtration par le poids de Bondarko
$$W_i{\bf H}_{\sigma}({\bf M} ,{\bf Z}(q))=\im({\bf H}_{\sigma}(w_{\leq i+2q}{\bf M} ,{\bf Z}(q))\rightarrow {\bf H}_{\sigma}({\bf M} ,{\bf Z}(q))).$$

Le principe (0.2) se traduit  en
 
\th PROPOSITION 3.3.4.
\enonce Pour chaque plongement  $\sigma$, et chaque entier positif $q$
les foncteurs  r\'ealisations de Betti ${\bf H}_{\sigma}(\ ,{\bf Z}(q))$ induisent des foncteurs
cohomologiques de la cat\'egorie des motifs g\'eom\'etriques sur $k$ vers la cat\'egorie des 
${\bf Z}$-modules filtr\'es de type fini. Si de plus le plongement 
est r\'eel, la r\'ealisation est munie d'une involution induite 
par la conjugaison complexe.
\endth

Comme dans le cas de la r\'ealisation de De Rham, nous avons
\th PROPOSITION 3.3.5. 
\enonce Soient $\M$ un motif g\'eom\'etrique effectif, $\M^{(\cdot)}$ un complexe des poids associ\'e \`a $\M$. La cohomologie du complexe total ${\bf K}_\sigma$ du bicomplexe $\r^\cdot\ho (t_\sigma(\M^{(\cdot)}),\Z(q))$ calcule la r\'ealisation de Betti de $\M$.
\endth

\goodbreak
\beginsection{4. R\'ealisation de Hodge}

\beginsection{4.1. Comparaison des r\'ealisations de De Rham et de Betti}

Fixons un plongement $\sigma : k \rightarrow \C$. Pour 
comparer $\ho_{D({\cal A}b)}(\M_{\sigma}(\C ), \Z (q))$ et 
$ \ho_{\dm(k)}(\M,  \om)$ pour tout motif $\M$,  nous partons des classes de 
Chern (2.1.5)  $\Z (q) \rightarrow \om$ et appliquons le foncteur 
de r\'ealisation topologique $t_{\sigma}$ (3.3) qui induit les fl\`eches
$$(2i\pi)^q \Z \rightarrow t_{\sigma} (\om)$$
dans la cat\'egorie des complexes de groupes ab\'eliens. En tensorisant avec le 
corps des nombres complexes, nous obtenons

\smallskip 
\th LEMME 4.1.1.
\enonce Le morphisme $\C \rightarrow t_{\sigma} (\om)$ est un quasi-isomorphisme 
de complexes de groupes ab\'eliens.
\endth
{\pc DEMONSTRATION}: au niveau des faisceaux, le foncteur de r\'ealisation 
topologique se 
factorise \`a travers le foncteur $t_{an}$ de (1.3.2.1). Le 
lemme de Poincar\'e implique que dans la cat\'egorie d\'eriv\'ee des faisceaux anaytiques, le complexe des diff\'erentielles holomorphes $\om_{an}$ est 
une r\'esolution du faisceau constant $\C$. Il induit un isomorphisme 
$\C\simeq  \om_{an}$  dans $\DMan$ qui se transporte dans $D(\widetilde{CW})$ 
puis $D({\cal A}b)$.
\hfill\cqfd

\goodbreak\smallskip
On en d\'eduit, pour tout complexe motivique $\M$,  une fl\`eche
$${\bf H}^{\,\cdot}_{DR} (\M) \otimes_k \C \rightarrow {\bf H}^{\,\cdot}_{\sigma} (\M,0) \otimes_{\Q} \C.$$
Par le th\'eor\`eme de De Rham, cette fl\`eche est un isomorphisme d'espaces 
vectoriels gradu\'es de dimension finie pour tout motif $\M(X)$ d'un sch\'ema 
lisse projectif $X$, isomorphisme compatible par construction avec l'action de 
la conjugaison complexe si le plongement $\sigma$ est r\'eel. Comme nous l'avons d\'ej\`a remarqu\'e, les filtrations par le poids co\"\i ncident (cf 3.3.2). Les motifs 
$\M(X)$ engendrant la cat\'egorie des motifs g\'eom\'etriques, on a 

\smallskip 
\th PROPOSITION 4.1.2.
\enonce Pour tout motif g\'eom\'etrique $\M$ de $\dg$, tout entier $q$ positif et tout entier $p$, on a un isomorphisme  de 
$\C$-espaces vectoriels filtr\'es (par la filtration par le poids) de dimension finie
$${\bf H}^p_{DR} (\M) \otimes_k \C \simeq {\bf H}^p_{\sigma} (\M,0) \otimes_{\Q} \C$$
compatible avec l'action de la conjugaison complexe si le plongement est r\'eel.
\endth

\beginsection{4.2. Th\'eorie de Hodge}

Une structure de Hodge mixte sur un sous-corps $k$ de $\C$ est la donn\'ee  d'un $\Z$-module de type fini $H_{\Z}$, d'une filtration croissante $W_n$ sur $H_{\bf Q}={\bf Q}\otimes_{\bf Z} H_{\Z}$, d'une filtration d\'ecroissante $F^p$ sur $H_k=k\otimes_{\bf Z}H_{\Z}$, tels que le syst\`eme des trois filtrations $(W, F ,\bar F)$ qui en d\'ecoule sur $H_\C={\bf C}\otimes_{\bf Z} H_{\Z}$ soit un syst\`eme de trois filtrations oppos\'ees [D71](1.2.13). La cat\'egorie des structures de Hodge mixtes est une cat\'egorie ab\'elienne (loc. cit. 2.3.5).

La filtration par le poids de Bondarko sur ${\bf H}^p_{\sigma}(\M, \Z)$, la filtration de Hodge sur ${\bf H}^p_{DR} (\M) $ et l'isomorphisme de comparaison De Rham-Betti permettent d'associer \`a un complexe motivique, objet de $\dg(k)$, une structure de Hodge mixte sur $k$.  Le th\'eor\`eme 0.1 se traduit en 

\th THEOREME 4.2.1.
\enonce Pour chaque plongement $\sigma: k \hookrightarrow \C$ , la donn\'ee du $\bf Z$-module ${\bf H}_{\sigma}(\M,\Z)$, de la filtration $W$ par le poids de Bondarko sur ${\bf Q}\otimes_{\bf Z}{\bf H}_{\sigma}(\M,\Z)$, du $k$-espace vectoriel filtr\'e par la filtration de Hodge $({\bf H}_{DR} (\M),F^\cdot)$, du complexifi\'e ${\bf H}_{\C}(\M) ={\bf C}\otimes{\bf H}_{DR} (\M)$, muni des trois filtrations $W_\cdot$, $F^\cdot$ et de sa complexe conjugu\'ee ${\bar F}^\cdot$, pour tout motif g\'eom\'etrique $\M$, induit un foncteur cohomologique de la cat\'egorie $\dg(k)$ des motifs
g\'eom\'etriques vers la cat\'egorie des structures de Hodge mixtes sur $k$.
\endth

 {\pc DEMONSTRATION}: 

Il suffit de traiter le cas o\`u $\M$ est un motif g\'eom\'etrique effectif. La d\'emonstration se fait en deux temps.

\th LEMME 4.2.2.
\enonce Si $\M$ est un motif de Chow, les donn\'ees du th\'eor\`eme fournissent une structure de Hodge pure.
\endth

{\pc DEMONSTRATION}: Si $\M$ est le motif d'un sch\'ema lisse projectif $X$, c'est un r\'esultat classique, puisque les r\'ealisations de $\M$ co\"\i ncident avec les r\'ealisations classiques de $X$.

Si $\M$ est une somme finie ou un facteur direct d'un $\M(X)$, pour $X$ projectif lisse, cela provient du th\'eor\`eme de Deligne que la cat\'egorie des structures de Hodge est ab\'elienne (cf. [D71]).

\th LEMME 4.2.3.
\enonce La donn\'ee des bicomplexes

{\leftskip1cm
\itemitem{$a)$} $K_B^{p,q}=\r^q\ho ({\bf M}_\sigma^{(-p)}({\bf C}),{\bf Z})$;   

\itemitem{$b)$} $K_{\bf Q}^{p,q}=\r^q\ho (\M_\sigma^{(-p)}({\bf C}),{\bf Q})$,  muni de la filtration par les poids;

\itemitem{$c)$} $K_{DR}^{p,q}=\r^q\ho(\M^{(-p)},{\bf \Omega})$,  muni de la premi\`ere filtration et de la filtration de Hodge,

\par}
est un $\bf Z$-complexe de Hodge mixte DG au sens de Deligne ([D74], 8.1.10).
\endth

{\pc DEMONSTRATION}: il suffit de v\'erifier que  pour chaque $p$, $K^{p,\cdot}$ est un complexe de Hodge mixte, ce qui provient du fait que chaque $M^{(-p)}$ est un motif de Chow et du lemme pr\'ec\'edent, qui permet de pr\'eciser que pour $p$ fix\'e, 
$K^{p,\cdot}$ est un complexe de Hodge pur de poids $0$.

\smallskip
Le th\'eor\`eme de Deligne (8.1.15 de [D74]) permet de conclure que pour chaque motif g\'eom\'etrique $\M$, les donn\'ees du th\'eor\`eme,  qui sont les cohomologies des donn\'ees du lemme,  d\'efinissent une structure de Hodge mixte.  On note  que le complexe total associ\'e au bicomplexe $K^{p,q}$, muni de sa filtration diagonale est bien le complexe ${\bf K}_{DR}$ muni de la filtration par le poids de Bondarko.
\smallskip
Sachant que la filtration par le poids est ind\'ependante des choix ([Bo10]), il nous faut v\'erifier 
\th PROPOSITION 4.2.4. 
\enonce
Si $\M$ est un motif g\'eom\'etrique, la filtration de Hodge  sur ${\bf H}_{DR}(\M)$ est ind\'ependante du choix du complexe des poids.
\endth
{\pc DEMONSTRATION}: Soient $\M^{(\cdot)}$ et  $\M'^{(\cdot)}$ deux complexes de poids du motif $\M$ que l'on peut supposer effectif. Bondarko a montr\'e en [Bo09] que les deux complexes \'etaient homotopes.  Les morphismes et les homotopies \'etant d\'efinis au niveau des motifs, ils induisent sur  ${\bf H}_{DR}(\M)$ des morphismes  compatibles aux filtrations de Hodge d\'efinies par troncature. Ainsi les deux complexes $\M^{(\cdot)}$ et  $\M'^{(\cdot)}$ d\'efinissent sur ${\bf H}_{\C}(\M)$ deux structures de Hodge mixtes et l'\'equivalence d'homotopie entre les deux complexes un automorphisme de ${\bf H}_{\C}(\M)$ compatible \`a toutes les filtrations: d'apr\`es [D71] cet automorphisme induit un isomorphisme de structures de Hodge mixtes.
\hfill\cqfd

\smallskip
On a ainsi construit un foncteur cohomologique ${\bf H}_{Hodge}$ de la cat\'egorie $\dge(k)$ dans la cat\'egorie des structures de Hodge mixtes. Pour conclure la d\'emonstration du th\'eor\`eme, il suffit de  pr\'eciser que le foncteur est multiplicatif et que les structures de Hodge des motifs de Tate sont inversibles, en rappelant
\smallskip 
\th{PROPOSITION} 4.2.5.
\enonce La r\'ealisation de Hodge des motifs de Tate $\Z (n)$ est la structure de Hodge pure de poids $2n$ suivante
$$\eqalign{ {\bf H}_{Hodge}^i({\bf Z}(n)) &=0 \quad\hbox{si}\  i\not =  0,\cr
{\bf H}_{Hodge}^0({\bf Z}(n)) &=({\bf H}^0_\sigma(\Z(n)),{\bf H}^0_{DR}(\Z(n)),\C\otimes{\bf H}^0_{DR}(\Z(n)),W_\cdot, F^\cdot,\bar F^\cdot)\cr
&=({\bf Z}, k, {\bf C}, W_\cdot, F^\cdot,\bar F^\cdot)}$$
avec pour filtration par le poids
$$\eqalign{ W_i {\bf H}^0_\sigma (\Z (n)) &=0\quad\hbox{ pour }i\leq 2n-1\cr 
 W_i {\bf H}^0_\sigma (\Z (n)) &= k \quad\hbox{ pour }i\geq 2n}$$
et pour filtration de Hodge 
$$\eqalign{ F^p {\bf H}^0_{DR}(\Z (n)) &=0\quad\hbox{ pour }p\geq n+1\cr 
 F^p {\bf H}^0_{DR} (\Z (n)) &= k \quad\hbox{ pour }p\leq n.}$$
\endth

{\pc DEMONSTRATION}: cela provient du triangle exact scind\'e [MVW Chap.15]
$$\M({\bf P}^{n-1}) \rightarrow \M({\bf P}^n) \rightarrow \Z (n)[2n] 
\rightarrow \M({\bf P}^{n-1})[1]$$ 
et du fait que l'espace projectif ${\bf P}^n$ a comme nombres 
de Hodge non nuls $h^{p,p}$ pour $0 \leq p\leq n$.
\hfill\cqfd

{\pc REMARQUE} 4.2.6. Remarquons que pour tout $q\in\bf Z$ et tout motif g\'eom\'etrique $\M$, 
$${\bf H}_\sigma(\M,{\bf Z}(q))={\bf H}_\sigma(\M(-q),{\bf Z})$$
et que cet isomorphisme est compatible avec la filtration par les poids d\'ecrite en 3.3.3. On peut ainsi d\'efinir une r\'ealisation de Hodge pour tout $q\in\bf Z$, en posant
$${\bf H}_{Hodge}(\M,q)={\bf H}_{Hodge}(\M(-q)).$$
\bigskip
\th{PROPOSITION} 4.2.7.
\enonce Si $\M= \M(X)$ est le motif d'un sch\'ema $X$, alors la structure de Hodge mixte d\'efinie en 4.2.1 co\"\i ncide avec celle d\'efinie par Deligne en [D74].
\endth

{\pc DEMONSTRATION} : si le sch\'ema $X$ est lisse et projectif, c'est vrai par construction (2.3.5). Si le sch\'ema $X$ est lisse, on le plonge dans un sch\'ema projectif lisse $\bar X$ tel que le sch\'ema compl\'ementaire $Y = {\bar X} - X$ soit un diviseur \`a croisements normaux \`a composantes $Y_n$ lisses, o\`u $n\leq N$. On choisit comme complexe des poids de $\M(X)$ le complexe d\'efini en (1.5.4), avec les notations de [D71]
$$\M(X)^{(i)} = \left\{\matrix{\M (\widetilde {Y^i})(i)[2i] & \hbox{ si }\   0\leq i \leq N\cr
0\hfill & \hbox{ sinon }\cr
}\right.$$
Les morphismes propres $i_m : \widetilde {Y^m} \rightarrow {\bar X}$ permettent de d\'efinir un bicomplexe de ${\cal O}_{\bar X}$-modules 
$i_{m\star}\Omega^p_{\widetilde Y^m}$ et un complexe  de ${\cal O}_{\bar X}$-alg\`ebres diff\'erentielles gradu\'ees
$$\oplus_{1\leq m \leq N} i_{m\star}\Omega_{\widetilde Y^m}^{\cdot\,-m}.$$

En [D71] Deligne d\'efinit la structure de Hodge mixte pour $X$ \` a partir du complexe de De Rham \`a p\^oles logarithmiques 
$\Omega^\cdot_{\bar X}\langle Y\rangle$ et construit une application (loc. cit. 3.1.5.2)
$$\Omega^\cdot_{\bar X}\langle Y\rangle \rightarrow \oplus_{1\leq m \leq N} i_{m\star}\Omega_{\widetilde Y^m}^{\cdot\,-m}(\epsilon^m)[-m]$$
qui induit un morphisme de bicomplexes
$$\eqalign{
\r\ho({\bar X}, \Omega^\cdot_{\bar X}\langle Y\rangle) \quad\longrightarrow\quad &\oplus_{1\leq m \leq N}\r\ho({\bar X}, i_{m\star}\Omega_{\widetilde Y^m}^\cdot (\epsilon^m)[-2m]\cr
&= \oplus_{1\leq m \leq N}\r\ho({\widetilde Y^m}, \Omega_{\widetilde Y^m}^\cdot(\epsilon^m)[-2m])\cr
&= \oplus_{1\leq m \leq N}\r^\cdot\ho(\M(X)^{(m)}, \om)\cr
}$$
compatible aux filtrations par le poids et de Hodge.

Ce morphisme  induit un isormophisme via le r\'esidu de Poincar\'e sur les termes $\gr_W$ et donc un isomorphisme entre les structures de Hodge respectives.

Si le sch\'ema $X$ est singulier, Deligne construit un sch\'ema simplicial lisse \`a partir d'\'eclatements, puis un $\Z$-complexe de Hodge mixte DG. On v\'erifie que par les triangles d'\'eclatements (1.3.5.4) on peut choisir un complexe des poids de telle sorte que le $\Z$-complexe de Hodge mixte DG construit en (4.2.3) soit quasi-isomorphe \`a celui de Deligne. On peut \'egalement utiliser le raisonnement de Huber (cf 4.4 ci-dessous).
\hfill\cqfd

\medskip
4.2.8. {\sl Application}

La th\'eorie  de Hodge  des sch\'emas projectifs lisses  complexes $X$ permet d'identifier $F^p{\Bbb H}^\cdot
(X, \Omega_X^\cdot)$ \`a ${\Bbb H}^\cdot(X, {}_{\sigma\geq p}\Omega_X)$ 
et le conoyau ${\Bbb H}^\cdot(X,\Omega^\cdot_X)/F^p{\Bbb H}^\cdot(X, \Omega^\cdot_X)$ \`a ${\Bbb H}^\cdot(X, \Omega_X^{<p})$. Cette propri\'et\'e s'\'etend aux motifs sur $\bf C$. Consid\'erons
le complexe $\om^{<p} $ de $D^-(\sh_{Nis}(\smc(\C)))$ 
$$ 0\rightarrow {\cal O}\rightarrow \om^1\rightarrow\cdots \rightarrow \om^{p-1}\rightarrow 0$$

\smallskip
\th {LEMME} 4.2.8.1.
\enonce Pour tout entier $p\geq 0$, pour  tout motif g\'eom\'etrique effectif de $\dme(\C)$ et tout complexe des poids $\M^{(\cdot)}$ de $\M$, l'hypercohomologie  du bicomplexe 
$$\r^\cdot\ho_{D(\sh_{Nis}(\smc(\C)))}(\M^{(\cdot)},\om^{<p}))$$
s'identifie au conoyau
$${\bf H}^\cdot_{DR}(\M)/F^p{\bf H}^\cdot_{DR}(\M).$$
\endth

{\pc DEMONSTRATION:} on plonge $\dme(\C)$ dans $D(\sh_{Nis}(\smc(\C)))$ o\`u l'on a, avec les notations \'evidentes, le triangle distingu\'e
$$ 0\rightarrow {}_{\sigma\geq p}\om \rightarrow \om \rightarrow \om^{<p} \rightarrow 0.\leqno (4.2.8.2)$$
On d\'eduit du cas des sch\'emas projectifs lisses que le triangle (4.2.8.2)  induit pour tout motif $\M$ de poids $0$ une suite exacte courte
$$0\rightarrow \ho_{D(\sh_{Nis}(\smc(\C)))}(\M, {}_{\sigma\geq p}\om)\rightarrow {\bf H}_{DR}(\M) \rightarrow \ho_{D^-(\sh_{Nis}(\smc(\C)))}(\M,\om^{<p})\rightarrow 0$$
avec $\ho_{D(\sh_{Nis}(\smc(\C)))}(\M,{}_{\sigma\geq p}\om) = F^p {\bf H}_{DR}(\M)$. Le lemme s'ensuit par une r\'ecurrence sur la longueur du complexe de poids de $\M$.\hfill\cqfd
\smallskip
{\pc REMARQUE}: le lemme signifie que bien que  $\om^{<p}$
ne soit pas un complexe motivique on peut donner un sens \`a $\r^\cdot\ho(\M, \om^{<p})$ en posant 
$$\r^\cdot\ho(\M, \om^{<p})={\rm Tot}\ \ \hbox{ c\^one}  (F^p\r^\cdot\ho_{\DM(\C)}(\M^{(\cdot)}, \om) \rightarrow \r^\cdot\ho_{\DM(\C)}(\M^{(\cdot)}, \om))$$
qui est bien d\'efini \`a homotopie pr\`es.

\medskip
4.2.9. Tous les r\'esultats concernant ${\bf H}_{\bf C}(\M)$ peuvent se voir dans $\dg(\C)$ via le changement de base, puis en g\'eom\'etrie analytique, via [GAGA], comme chez Deligne. En particulier, le lemme 4.2.8.1. se r\'einterpr\`ete analytiquement en
$$ {\Bbb H}^n(\r^\cdot\ho_{D^-_{an}}(t_{an}({\widetilde\M^{(\cdot)}}), \om_{an}^{<q}))= {\bf H}^n_{DR}(\M)/F^q{\bf H}^n_{DR}(\M).$$

\beginsection{4.3. R\'ealisation  de Deligne-Beilinson.}

Les r\'esultats pr\'ec\'edents permettent de g\'en\'eraliser aux motifs g\'eom\'etriques la d\'efinition de la cohomologie de Deligne-Beilinson. Pour les d\'efinitions des cohomologies de Deligne et Deligne-Beilinson des sch\'emas, nous nous r\'ef\'erons \`a [Sc88] et [EV88].

Rappelons que pour un sch\'ema $X$ lisse sur $\C$, le complexe de Deligne est d\'efini pour tout entier $q$ comme le complexe de faisceaux en topologie analytique $D_X(q)$
$$0\rightarrow (2i\pi)^q \Z \rightarrow {\cal O}_X^h \rightarrow \Omega_X^{1,h} \dots \rightarrow \Omega_X^{q-1,h} \rightarrow 0 $$
o\`u le terme $(2i\pi)^q \Z = \Z(q)_{an}$ est en degr\'e $0$.
Ce complexe s'\'ecrit 
$$D_X(q) = \ \hbox{ c\^one} \  (\Z(q)_{an} \rightarrow \Omega_X^{\cdot,h}/F^q\Omega_X^{\cdot,h}) [-1]$$
et on pose $H^\cdot_D (X,q) = H_{an}^\cdot (X,D_X(q))$.

Si le sch\'ema $X$ est projectif, la th\'eorie de Hodge  et [GAGA] permettent d'inclure $H^\cdot_D (X,q)$ dans une longue suite exacte
$$\dots \rightarrow H_D^n(X,\Z(q)) \rightarrow H^n(X(\C), \Z(q))\rightarrow H^n_{DR} (X(\C))/F^q \rightarrow H_D^{n+1}(X, \Z(q))\dots $$
o\`u l'on voit que deux termes sur trois ont une interpr\'etation motivique. Si le sch\'ema n'est pas projectif, Beilinson [B84] a  d\'ecrit le complexe en terme de complexe de De Rham \`a p\^oles logarithmiques pour d\'efinir une cohomologie de $X$ qui s'inscrit dans la m\^eme suite exacte longue. Comme dans le cas de la r\'ealisation de Hodge on va r\'esoudre un motif par son complexe des poids.

Dans la cat\'egorie d\'eriv\'ee des faisceaux analytiques sur $\smc(\C)$, notons  $\om_{an}^{<p}$ le complexe 
$$0\rightarrow {\cal O}^h \rightarrow \om^{1,h} \dots \rightarrow \om^{p-1,h} \rightarrow 0 $$
o\`u l'anneau des fonctions holomorphes est en degr\'e $0$ et posons 
$$D(q) = \ \hbox{ c\^one} \ \ (\ \Z(q)_{an} \rightarrow \om_{an}^{<q})[-1].$$

\smallskip 
Soit $\M$ un motif g\'eom\'etrique effectif et $\M^{(\cdot)}$ un complexe des poids de $\M$. Par le th\'eor\`eme d'invariance d'homotopie, on identifie ces motifs de $\dm(\C)$ \`a des complexes $\widetilde\M$ et   $\widetilde{\M^{(\cdot)}}$ de $D^-(\sh_{Nis}(\smc(\C)))$.  Consid\'erons le complexe 
$$DB_q(\M) = {\rm Tot} \ \ \r^\cdot\ho_{D^-_{an}}(t_{an} (\widetilde{\M^{(\cdot)}}), D(q)).$$ 

\th PROPOSITION 4.3.3.
\enonce Pour tout motif g\'eom\'etrique effectif $\M$ le complexe $DB_q(\M)$ est bien d\'efini \`a homotopie pr\`es et sa cohomologie not\'ee $H_D(\M,q)$ est isomorphe quand $\M$ est le motif d'un sch\'ema lisse complexe $X$ \`a la cohomologie de Deligne-Beilinson  de $X$.
\endth

{\pc DEMONSTRATION}:  le complexe $DB_q(\M)$ est le c\^one de deux complexes bien d\'efinis \`a homotopie pr\`es: un calcule la r\'ealisation de Betti (3.1.7) et le deuxi\`eme a \'et\'e consid\'er\'e en 4.2.9. Par 3.3.5, on a des isomorphismes
$${\bf H}^n(\M(\C),\Z(q))=\ho_{D^-_{an}}(t_{an}({\widetilde\M}), \Z(q)_{an} [n])={\Bbb H}^n\r^\cdot\ho_{D^-_{an}}(t_{an}({\widetilde\M^{(\cdot)}}), \Z(q)_{an}))$$  et la suite exacte longue 
$$\dots \rightarrow H^n(DB_q(\M)) \rightarrow {\bf H}^n(\M(\C),\Z(q))\rightarrow{\Bbb H}^n \r^\cdot\ho_{D^-_{an}}(t_{an}({\widetilde\M^{(\cdot)}}), \om_{an}^{<q}) \rightarrow H^{n+1}(DB_q(\M)) \rightarrow\dots$$
qui, pour $\M$  motif d'un sch\'ema lisse $X$,  s'identifie \`a celle d\'efinissant la cohomologie de Deligne-Beilinson.\hfill\cqfd
\smallskip

Le produit sur le complexe de Deligne $ D(q) \otimes D(q') \rightarrow D(q+q')$ est d\'efini \`a partir du produit des diff\'erentielles et commute aux transferts: 
$$x\cup y =\left\{ \matrix{
x.y &{\rm  si}\ &i=0 \cr
x\wedge dy \hfill &{\rm si} &\ i>0 \ {\rm et}\  j=q'\cr
0&{\rm sinon}\cr
}\right.$$
avec $x\in D(q)^i$ et $y \in D(q')^j$.
\smallskip
Le foncteur cohomologique $H_D$ est fonctoriel et induit un foncteur
de la cat\'egorie des motifs g\'eom\'etriques sur $\C$ vers la cat\'egorie des anneaux gradu\'es. 
\smallskip
Si $\sigma : k\hookrightarrow \C$ est un plongement, on appelle r\'ealisation de Deligne-Beilinson le foncteur cohomologique $H_{D,\sigma} = H_D \circ t_{\sigma}$
de $\dg(k)$ vers la cat\'egorie des anneaux gradu\'es.
Si le plongement $\sigma$ est r\'eel, on montre comme en (3.2) que la r\'ealisation  $H_{D,\sigma}(\M)$ est muni d'une action de la conjugaison complexe $F_{\infty}$. On appelle r\'ealisation de Deligne-Beilinson r\'eelle le foncteur des invariants par $F_{\infty}$ de $H_{D,\sigma}$. 

\bigskip
{\pc REMARQUE}:  dans la cat\'egorie des complexes motiviques sur $\C$, nous avons d\'efini en (2.1.5) des classes de Chern $\Z(q) \rightarrow \om$ qui se factorisent \`a travers $\Z(q) \rightarrow{}_{\tau{\leq q}}\om$; la projection $\Z(q) \rightarrow\om^{<q}$ est triviale si l'on reste dans une cat\'egorie de faisceaux alg\'ebriques, d'o\`u la n\'ecessit\'e de passer par les r\'ealisations analytiques.
En effet, pour $q=1$, il n'existe pas de morphisme alg\'ebrique non trivial entre le groupe multiplicatif et le groupe additif. Pour d\'efinir une r\'ealisation \`a la Deligne-Beilinson, il faut pouvoir construire une application logarithme. 

\beginsection{4.4. Comparaison avec les constructions de Huber}

Annette Huber construit ses r\'ealisations en \'etendant \`a la cat\'egorie 
$\dg(k)$ des motifs g\'eom\'etriques les foncteurs additifs qu'elle a 
d\'efinis en [H95]
$$\tilde R :\sm (k) \rightarrow C^+ ({\cal A})$$
de la cat\'egorie des sch\'emas lisses sur $k$ dans la cat\'egorie des complexes
born\'es inf\'erieurement d'une cat\'egorie ab\'elienne $\Q$-lin\'eaire 
${\cal A}$ (cf [H00], [H04] (theorem B.2.2)). 
Il est clair que, rationnellement,  nos r\'ealisations co\"\i ncident avec les siennes pour les motifs des vari\'et\'es lisses et il reste \`a
 v\'erifier que nos constructions s'\'etendent de la m\^eme fa\c con que la sienne. Comme elle le pr\'ecise, il est important que 
les 
foncteurs $\tilde R$ soient \`a valeurs dans une cat\'egorie de complexes et non dans
une cat\'egorie homotopique ou d\'eriv\'ee. Elle-m\^eme a construit, pour tout sch\'ema lisse,  des complexes [H95]
$\tilde R_{sing} (X)$ et $\tilde R_{dR} (X)$.

Fixons-nous des r\'esolutions 
injectives respectives $J_{DR}^{\cdot}$ du complexe de De Rham 
$\om$ et 
d\'efinissons, pour un plongement $\sigma \hookrightarrow \C$ fix\'e,
 les foncteurs
 $$\matrix{
 R_{DR} & : & \sm (k)  & \rightarrow  & \C^+ ({\cal F}_k) \cr
  & & X & \mapsto & \ho_{C(\sh _{Nis}(\smc k))} (C_*(\Ztr (X)),  
J_{DR}^{\cdot})\cr
  }$$
  o\`u ${\cal F}_k$ est la cat\'egorie des $k$-espaces vectoriels filtr\'es,
$$\matrix{
  R_B& : & \sm k& \rightarrow  & C^+({\cal E}_{\Q})\cr
  & & X & \mapsto &  \ho_{C({\cal A}b)}( {\rm Map} (\Delta_{\rm top}^{\cdot}, 
\coprod_{d\geq 0} S^d X(\C))^+, \Q),
  }$$
   o\`u ${\cal E}_{\Q}$ est la cat\'egorie des $\Q$-espaces vectoriels.
   
   Pour tout sch\'ema lisse $X$ le dual du morphisme de Dold-Thom fournit un quasi-isomorphisme
$$u_B(X): \tilde R_{sing} (X) \simeq  R_B (X),$$ le complexe $\tilde R_{dR} (X)$ construit \`a partir des diff\'erentielles \`a 
p\^oles logarithmiques  est une r\'esolution du motif de De Rham et induit un quasi-isomorphisme
$$R_{DR} (X)\simeq\tilde R_{dR} (X)$$
compatible aux filtrations par les poids (2.2) et de Hodge (2.3).

L'essence des constructions de Huber ([H00](2,3)) est d'\'etendre les foncteurs
 $\tilde R$ 
de $\sm k$ \`a la cat\'egorie $\smc k$, la proposition (loc. cit. Prop.2.1.2) 
assurant
qu'une fois choisie une r\`egle de signes dans les multicomplexes, un foncteur
de $\smc k$ dans $C^+({\cal A})$ v\'erifiant les bonnes propri\'et\'es 
s'\'etend naturellement en un foncteur $\dg(k) \rightarrow 
D^+({\cal A})$.

Par construction, nos foncteurs $R$ s'\'etendent en des foncteurs sur $\smc (k)$.
 Il nous faut montrer que ces derniers co\"\i ncident avec les foncteurs 
 \'etendus par Huber ([H00], theorem 2.1.6 et [H04] 
B.2.2), ce qui se ram\`ene \`a v\'erifier que nous avons d\'efini les m\^emes transferts.

\smallskip 
\th LEMME 4.4.1.
\enonce Pour toute correspondance $\alpha$ entre les sch\'emas lisses $X$ et $Y$
, on a un diagramme commutatif
$$\matrix{
\tilde R_{sing} (Y) & \hfl{\alpha^*_{sing}}{} & \tilde R_{sing} (X) \cr
\vfl{}{u_B(Y)} && \vfl{}{u_B(X)}\cr
R_B(Y)& \hfl{\alpha^*_B}{}& R_B(X)\cr
}$$
\endth
{\pc DEMONSTRATION}: la construction des transferts sur $\tilde R_{sing}$ est 
explicit\'ee dans la d\'emonstration des th\'eor\`emes 2.1.3 et 2.1.6 de [H95].
Celle des transferts sur $R_B$ est, avant application du foncteur de 
r\'ealisation 
topologique, issue de l'article [SV96] de Suslin et 
Voevodsky.
Dans les deux cas on se ram\`ene \`a une correspondance \'el\'ementaire  de 
support  qu'on peut supposer normal, et m\^eme un recouvrement 
g\'en\'eriquement galoisien chez Huber ou pseudo-galoisien chez 
Suslin et Voevodsky (les deux notions sont \'equivalentes en caract\'eristique $0$). La m\^eme trace permet alors de d\'efinir les 
transferts qui co\"\i ncident.
\hfill\cqfd

La  d\'emonstration se transpose au cas des r\'ealisations de De Rham et de Hodge puisque les transferts ([LW09]), puis les filtrations, ont \'et\'e construits au niveau des complexes. La naturalit\'e des isomorphismes de comparaison permet de conclure

\smallskip 
\th PROPOSITION 4.4.2.
\enonce Nos foncteurs de r\'ealisations ${\bf H}_{DR}$ et ${\bf H}_{\sigma}$ 
restreints 
\`a la cat\'egorie $\dg(k)$ des motifs g\'eom\'etriques co\"\i ncident
respectivement avec les composantes de Rham et singuli\`ere du 
foncteur de r\'ealisation mixte de Huber [H00].
\endth

Par ailleurs, comme l'a pr\'ecis\'e Bondarko ([Bo09], 7.4), il n'existe qu'une seule filtration par le poids sur une r\'ealisation de $\dg$ vers une cat\'egorie \`a coefficients rationnels. Nos filtrations par le poids co\"\i ncident donc avec celles de [H00], de m\^eme que les r\'ealisations de Hodge.

\biblio {Bibliographie}
\ref{SGA4} \AUTHOR = { Artin, M. - Grothendieck, A. - Verdier, J. L} 
\TITLE = {Th\'eorie des topos et cohomologie \'etale des sch\'emas.
     S\'eminaire de G\'eom\'etrie Alg\'ebrique du Bois-Marie
              1963--1964 (SGA 4)         {T}ome 1: {T}h\'eorie des topos}\SERIES 
= {Lecture Notes in Mathematics}\VOLUME= {269}
 \PUBLISHER = {Springer-Verlag}\ADDRESS = {Berlin}\YEAR = {1972}\PAGES = 
{xix+525}\book

\ref{A10}
\AUTHOR = {Ayoub, Joseph}
     \TITLE = {Note sur les op\'erations de {G}rothendieck et la
              r\'ealisation de {B}etti}
   \JOURNAL = {J. Inst. Math. Jussieu}
    \VOLUME = {9}
      \YEAR = {2010}
    \NUMBER = {2}
     \PAGES = {225--263}\article

\ref{B84}
    \AUTHOR = {Be{\u\i}linson, A. A.}
     \TITLE = {Higher regulators and values of {$L$}-functions}
 \BOOKTITLE = {Current problems in mathematics, {V}ol. 24}
    \SERIES = {Itogi Nauki i Tekhniki}
     \PAGES = {181--238}
 \PUBLISHER = {Akad. Nauk SSSR Vsesoyuz. Inst. Nauchn. i Tekhn. Inform.}
   \ADDRESS = {Moscow}
      \YEAR = {1984}\incollection

\ref{BV08}\AUTHOR = {
Be{\u\i}linson, Alexander - Vologodsky, Vadim}\TITLE = {A DG guide to Voevodsky's motives}
\JOURNAL = {Geom. Funct. Anal.}\VOLUME = {17} 
\YEAR = {2008}\NUMBER={6}\PAGES = { 1709--1787}\article

\ref{Bo09}\AUTHOR = {Bondarko, Mikhael V}\TITLE = {Differential graded motives: weight complex, weight filtrations and spectral sequences for realizations; Voevodsky versus Hanamura }\JOURNAL = { J. Inst. Math. Jussieu   }\VOLUME = { 8} 
\YEAR = {2009}\NUMBER={1}\PAGES = {39--97}\article

\ref{Bo10}\AUTHOR = {Bondarko, Mikhael V}\TITLE = {Weight structures, weight filtrations, weight spectral sequences, and weight complexes (for motives and spectra)  }\JOURNAL = {J. $K$-Theory} \VOLUME ={6} \YEAR = {2010} \NUMBER = {3} \PAGES = {387--504} \article

\ref{CE} \AUTHOR = {Cartan, Henri - Eilenberg, Samuel}\TITLE = {Homological algebra}\SERIES = {Princeton Landmarks in Mathematics}
      \PUBLISHER = {Princeton University Press}\ADDRESS = {Princeton, NJ}\YEAR = {1999}\book

\ref{CD09}\AUTHOR = {Cisinski, Denis-Charles - D\'eglise, Fr\'ed\'eric} \TITLE 
=  {Local and stable homological algebra in Grothendieck abelian
              categories} \JOURNAL = {Homology, Homotopy and Applications}
\VOLUME = {11}
\YEAR = {2009}
    \NUMBER = {1}
\PAGES = {219--260}
\article 

\ref{D71}
\AUTHOR = {Deligne, Pierre}
     \TITLE = {Th\'eorie de Hodge II}
    \JOURNAL= { Inst. Hautes \'etudes Sci. Publ. Math.}
    \VOLUME = {40}\YEAR = {1971}\PAGES = {5--57}\article

\ref{D74}
\AUTHOR = {Deligne, Pierre}
     \TITLE = {Th\'eorie de Hodge III}
    \JOURNAL= { Inst. Hautes \'etudes Sci. Publ. Math.}
    \VOLUME = {44}\YEAR = {1974}\PAGES = {5--77}\article

\ref{EV88}
\AUTHOR = {Esnault, H{\'e}l{\`e}ne and Viehweg, Eckart} \TITLE={Deligne-Be{\u\i}linson cohomology}
\BOOKTITLE = {Be\u\i linson's conjectures on special values of {$L$}-functions}
\SERIES = {Perspect. Math.}
\VOLUME = {4}
\PAGES = {43--91}
\PUBLISHER = {Kluwer Acad. Publ.}\ADDRESS = {Boston, MA}\YEAR = {1988}\incollection

\ref{FV}\AUTHOR = {Friedlander, Eric M. - Voevodsky, Vladimir}\TITLE = {Bivariant cycle cohomology}\BOOKTITLE = {Cycles, transfers, and motivic homology 
theories}\SERIES = {Ann. of Math. Stud.}\VOLUME = {143}\PAGES = 
{138--187}\PUBLISHER = {Princeton Univ. Press}\ADDRESS = {Princeton, NJ}\YEAR = 
{2000}\incollection

\ref{Go58}\AUTHOR = {Godement, Roger}\TITLE = {Topologie alg\'ebrique et 
th\'eorie des faisceaux}
 \SERIES = {Actualit\'es scientifiques et industrielles 1252. Publications de 
 l'Institut de Math\'ematique de l'Universit\'e de Strasbourg} \VOLUME={XIII}
 \PUBLISHER ={Hermann} \ADDRESS= {Paris}  \YEAR = {1858}\book

\ref{EGA IV}\AUTHOR = {Grothendieck, A.}
     \TITLE = {El\'ements de g\'eom\'etrie alg\'ebrique. {IV}. \'{E}tude
              locale des sch\'emas et des morphismes de sch\'emas {IV}}
  \JOURNAL = {Institut des Hautes \'Etudes Scientifiques. Publications
              Math\'ematiques}
    \NUMBER = {32}
      \YEAR = {1967}
     \PAGES = {361}\article

\ref{H95} \AUTHOR = {Huber, Annette} 
\TITLE = {Mixed motives and their realization in derived categories} 
\SERIES = {Lecture Notes in Mathematics} \VOLUME = {1604}
\PUBLISHER = {Springer-Verlag}\ADDRESS = {Berlin}\YEAR = {1995}\book

\ref{H00} \AUTHOR = {Huber, Annette}\TITLE = {Realization of {V}oevodsky's 
motives} \JOURNAL = {Journal of Algebraic Geometry}\VOLUME = {9}\YEAR = {2000} 
\NUMBER = {4}\PAGES = {755--799}\article	

\ref{H04}
\AUTHOR = {Huber, A.}
     \TITLE = {Corrigendum to: ``{R}ealization of {V}oevodsky's motives''
              [{J}. {A}lgebraic {G}eom. {\bf 9} (2000), no. 4, 755--799]}\JOURNAL  = {Journal of Algebraic Geometry}\VOLUME = 
{13}\YEAR = {2004}\NUMBER = {1}\PAGES = {195--207}\article

\ref{J90} \AUTHOR = {Jannsen, Uwe}
     \TITLE = {Mixed motives and algebraic {$K$}-theory}
    \SERIES = {Lecture Notes in Mathematics}
    \VOLUME = {1400}
 \PUBLISHER = {Springer-Verlag}
   \ADDRESS = {Berlin}
      \YEAR = {1990}\book

\ref{K73}\AUTHOR = {Knighten, Carol M.}
\TITLE = {Differentials on quotients of algebraic varieties}
\JOURNAL = {Transactions of the American Mathematical Society}
\VOLUME = {177}
\YEAR = {1973}\PAGES = {65--89}\article 

\ref{L08} \AUTHOR= {Lecomte, Florence} \TITLE= {R\'ealisation de Betti des motifs de 
Voevodsky}\JOURNAL = {Comptes rendus - Mathematique}\VOLUME = {346}
\YEAR = {2008}\PAGES = {1083--1086}\article
  
\ref{LW09} \AUTHOR={Lecomte, Florence - Wach, Nathalie}\TITLE = {Le complexe 
motivique de De Rham}\JOURNAL = {Manuscripta Math.}\VOLUME = {129}
\YEAR = {2009}\PAGES = {75--90}\article 

\ref{Ma68}\AUTHOR = {Manin, Ju. I}\TITLE = {Correspondences, motifs and 
monoidal transformations. 
(Russe)}\JOURNAL = { Izv. Akad. Nauk SSSR Ser. Mat. }\VOLUME = {32}\YEAR= 
{1968}
\PAGES= {1223--1244}\article

\ref{MVW}\AUTHOR = {Mazza, Carlo - Voevodsky, Vladimir - Weibel, Charles}\TITLE 
= {Lecture notes on motivic cohomology}\SERIES = {Clay Mathematics 
Monographs}\VOLUME = {2}\PUBLISHER = {American Mathematical Society}\ADDRESS = 
{Providence, RI}\YEAR = {2006}\book

\ref{N89}\AUTHOR = {Nisnevich, Ye. A}
\TITLE = {The completely decomposed topology on schemes and associated descent 
spectral sequences in algebraic {$K$}-theory}\BOOKTITLE = {Algebraic $K$-theory: 
connections with geometry and topology (Lake Louise, AB, 1987)}\SERIES = {NATO 
Adv. Sci. Inst. Ser. C Math. Phys. Sci.}\VOLUME = {279}\PAGES = 
{241--342}\PUBLISHER = {Kluwer Acad. Publ.}\ADDRESS = {Dordrecht}\YEAR = 
{1989}\incollection

\ref{Sc88}\AUTHOR = {Schneider, Peter}
\TITLE = {Introduction to the {B}e\u\i linson conjectures}
\BOOKTITLE = {Be\u\i linson's conjectures on special values of {$L$}-functions}
\SERIES = {Perspect. Math.}
\VOLUME = {4}
\PAGES = {1--35}
\PUBLISHER = {Kluwer Acad. Publ.}\ADDRESS = {Boston, MA}\YEAR = {1988}\incollection

\ref{GAGA} \AUTHOR = {Serre, Jean-Pierre}\TITLE = {G\'eom\'etrie alg\'ebrique et 
g\'eom\'etrie analytique}\JOURNAL = {Universit\'e de Grenoble. Annales de 
l'Institut Fourier}\VOLUME = {6}\YEAR = {1955--1956}\PAGES = {1--42}\article

\ref{SV96} \AUTHOR = {Suslin, Andrei - Voevodsky, Vladimir}
\TITLE = {Singular homology of abstract algebraic varieties}
\JOURNAL = {Inventiones Mathematicae}\VOLUME = {123}\YEAR = {1996}
\NUMBER = {1}\PAGES = {61--94}\article

\ref{Ver77} \AUTHOR = {Verdier, Jean-Louis} \TITLE = {Cat\'egories d\'eriv\'ees, 
\'etat 0}\BOOKTITLE = {Cohomologie \'etale, S\'eminaire de G\'eom\'etrie 
Alg\'ebrique du Bois-Marie SGA 4${1\over 2}$} \SERIES = {Lecture Notes in 
Mathematics} \VOLUME = {569} \PAGES = {262--311}\PUBLISHER = 
{Springer-Verlag}\ADDRESS = {Berlin}
      \YEAR = {1977}\incollection

\ref{Ver96}
\AUTHOR = {Verdier, Jean-Louis} \TITLE = {Des cat\'egories d\'eriv\'ees des cat\'egories ab\'eliennes} \JOURNAL = {Ast\'erisque} \NUMBER = {239} \YEAR = {1996} \PAGES = {xii+253 pp. (1997)}\article

\ref{V-TCM}\AUTHOR = {Voevodsky, Vladimir}\TITLE = {Triangulated categories of 
motives over a field}\BOOKTITLE = {Cycles, transfers, and motivic homology 
theories}\SERIES = {Ann. of Math. Stud.}\VOLUME = {143}\PAGES = 
{188--238}\PUBLISHER = {Princeton Univ. Press}\ADDRESS = {Princeton, NJ}\YEAR = 
{2000}\incollection

\bye